\newif\ifpomstemplate
\IfFileExists{poms1_V1.cls}{\pomstemplatetrue}{\pomstemplatefalse}
\ifpomstemplate
  \documentclass[poms,final]{poms1_V1}
  \OneAndAHalfSpacedXI
\else
  \documentclass[11pt]{article}
  \usepackage[margin=1in]{geometry}
  \usepackage{setspace}
\fi

\usepackage[T1]{fontenc}
\usepackage[utf8]{inputenc}
\usepackage{booktabs}
\usepackage{longtable}
\usepackage{enumitem}
\usepackage{float}
\usepackage{tabularx}
\usepackage{array}
\usepackage{graphicx}
\usepackage{caption}
\usepackage{natbib}
\usepackage{url}
\usepackage{xcolor}
\usepackage{tikz}
\usetikzlibrary{arrows.meta,positioning,fit,shapes.geometric}
\makeatletter
\@ifpackageloaded{hyperref}{}{\usepackage[hidelinks]{hyperref}}
\makeatother
\bibpunct[, ]{(}{)}{,}{a}{}{,}
\def\bibfont{\small}
\def\bibsep{\smallskipamount}
\newcolumntype{Y}{>{\raggedright\arraybackslash}X}
\newcolumntype{L}[1]{>{\raggedright\arraybackslash}p{#1}}
\newcolumntype{P}[1]{>{\raggedright\arraybackslash}p{#1}}
\tikzset{
  process/.style={draw, rounded corners, align=center, inner sep=5pt,
    minimum height=9mm, fill=blue!5},
  gate/.style={draw, diamond, aspect=1.9, align=center, inner sep=2pt,
    fill=orange!10},
  flow/.style={-{Latex[length=2mm]}, thick},
  loopflow/.style={-{Latex[length=2mm]}, thick, rounded corners}
}

\ifpomstemplate
  \renewcommand{\theARTICLETOPRIGHT}{}
\fi

\begin{document}

\ifpomstemplate
  \RUNAUTHOR{Liu, Lin, Chen, Wu, and Shen}
  \RUNTITLE{Machine Learning for Scheduling Decision Systems}
  \TITLE{Machine Learning for Scheduling Decision Systems: A Critical Review of Architecture, Assurance, and Deployment}
  \ARTICLEAUTHORS{
  \AUTHOR{
    Anbang Liu\textsuperscript{1},
    Shaochong Lin\textsuperscript{1},
    Jingchuan Chen\textsuperscript{1,2},
    Peng Wu\textsuperscript{3},
    and Zuo-Jun Max Shen\textsuperscript{4,5}
  }

  \AFF{
    \textsuperscript{1}
    Department of Data and Systems Engineering,
    University of Hong Kong, Hong Kong, China
    (\EMAIL{anbang@hku.hk}; \EMAIL{shaoclin@hku.hk})
    \\[2pt]
    \textsuperscript{2}
    Business School, Central University of Finance and Economics,
    Beijing, China
    (\EMAIL{jingchuan.chen@cufe.edu.cn})
    \\[2pt]
    \textsuperscript{3}
    School of Business, Sichuan University,
    24 South Section 1, 1st Ring Road,
    Chengdu, China
    (\EMAIL{pengwu@scu.edu.cn})
    \\[2pt]
    \textsuperscript{4}
    Faculty of Engineering and Faculty of Business and Economics,
    University of Hong Kong, Hong Kong, China
    (\EMAIL{maxshen@hku.hk})
    \\[2pt]
    \textsuperscript{5}
    College of Engineering, University of California, Berkeley,
    Berkeley, CA, USA
  }
}

\ABSTRACT{
Machine learning can support scheduling through prediction, search guidance, or schedule formation, but model-level evaluations do not reveal the downstream work, technical authority, and controls required to release a decision. We conduct a critical integrative review combining structured candidate identification with purposive full-text synthesis. The unit of analysis is the complete reported scheduling decision system---from problem specification through release and conditional recovery.

We develop a system-level taxonomy that distinguishes what a learned output does, how schedule content is formed, and which components exercise binding control along the normal release path. Learning and adaptation are classified separately from decision assurance, while organizational decision rights remain distinct from technical authority. Applying the taxonomy shows that reported timing results and solver guarantees depend on downstream work and the decision space left open after learned commitments. It also distinguishes transfer of learned capability from architectural reuse and operational maturity from automated release.

The evidence supports selected quality--computation-time trade-offs, bounded transfer of retained learned capability, and performance under specified operating regimes; it also documents operating use in several configurations. It does not provide a system-equivalent ranking of learning and optimization, evidence of general cross-task transfer, or a common basis for lifecycle economics and long-run field performance.

These findings motivate four operations-management propositions linking relative lifecycle value to effective reuse, the fit between the allocation of technical authority and the form of change, deadline-feasible assurance and recovery, and the alignment of organizational rights with information and accountability. Solver-led configurations, configurations with shared technical authority, and model-led configurations are alternative designs rather than maturity stages; greater model performance alone does not justify greater release authority.
}
  

  \maketitle
\else
  \title{Machine Learning for Scheduling Decision Systems: A Critical Review of Architecture, Assurance, and Deployment}
  \author{Anbang Liu, Shaochong Lin, Jingchuan Chen, Peng Wu, and Max Z.J. Shen}
  \date{}
  \maketitle
\fi

\section{Introduction}
\label{sec:introduction}

Scheduling allocates scarce capacity over time to meet operational commitments.
Temporal and resource interdependence make many stylized scheduling problems
NP-hard, while operating settings add changing information and short decision
windows. Classical research provides exact and approximate methods
\citep{pinedo2012scheduling}. In recurring operations, however, organizations
must decide not only how to compute a schedule, but also how schedule formation,
checking, release, and recovery should be divided among learned models,
optimizers, and responsible personnel. Learning-augmented scheduling is
therefore an operations-management system-design problem, not only an algorithm
choice.

Learning in scheduling has a long lineage, from reinforcement-learning
dispatching and shop-control policies to neural combinatorial optimization (NCO)
\citep{zhang1995rljobshop,riedmiller1999neural,aydin2000dynamic,
gabel2012distributed,vinyals2015pointer,kool2019attention,kwon2020pomo,
kwon2021matnet}. At inference, learning may supply information to a retained
optimizer, propose schedule content, or form and refine a schedule. The
capability supporting these roles may be created from reference targets or
direct evaluation
\citep{li2024imitation,corsini2024selflabeling,kotary2022fast,
zhang2020dispatch}. Two studies can therefore report similar model quality or
inference time while leaving different downstream search, checking, and release
work. Model-level results need not imply system-equivalent performance, cost,
technical authority, or assurance.

Existing reviews organize learning for combinatorial and constrained optimization or examine particular scheduling models and operating contexts
\citep{bengio2021methodological,kotary2021survey,cappart2023combinatorial, smit2025gnnsurvey,lv2025drlreview,khadivi2025drlreview,zhang2025rlreview, goppert2024online,zhang2024resilientreview,stockermann2025deploymentreview}. Taken together, these reviews leave unresolved how learned and retained components are assembled into a complete scheduling decision system. We conduct a critical integrative review that takes the complete reported system---from problem specification through release and conditional recovery---as its conceptual unit. We ask how these systems allocate work and binding technical authority, what the evidence supports at that boundary, and when alternative configurations create lifecycle value. Our central argument is that learning reallocates scheduling work and control rather than uniformly replacing optimization. Table~EC.2 in the electronic companion (EC) positions
this system-level view relative to prior reviews.

Our first contribution is a layered taxonomy of learning-augmented scheduling
systems. Within inference architecture, it separates learning interface,
solution formation, and normal-path technical authority. Learning and adaptation, decision assurance, and organizational decision rights are recorded separately, with conditional recovery distinguished from the normal release
path. This structure compares systems that use similar models but allocate
downstream work and binding control differently.

Our second contribution is a claim-level evidence audit at the same system
boundary. Applying the taxonomy shows why inference time is not complete-system
response time, why solver guarantees cover only the decision space left open
after learned commitments, why architectural reuse is not transfer of learned
capability, and why operating use does not by itself establish automated
release. The evidence supports selected quality--computation-time trade-offs,
bounded transfer of retained learned capability, performance under specified
operating regimes, and the operational possibility of several configurations.
It does not provide a system-equivalent ranking of learning and optimization or
a common basis for general cross-task transfer, lifecycle economics, and
long-run field performance.

Our third contribution is an operations-management contingency framework for
investment, allocation of technical authority, automated release, and
organizational governance. The audit identifies four unresolved relationships,
while operations-management reasoning and relevant theory supply their proposed
direction. The resulting propositions treat solver-led configurations,
configurations with shared technical authority, and model-led configurations as
alternative designs whose value depends on the operating setting.

\subsection{Review Scope and Evidence Boundaries}
\label{subsec:review-scope}

The review combines structured literature identification with a purposively
assembled full-text analytical library. A 46-query search across Semantic
Scholar, DBLP, arXiv, OpenReview, and PubMed, with an August 26, 2026 cutoff,
produced 2,408 consolidated candidate substantive-work families. This search
supports discovery rather than a prevalence estimate. The 148-work analytical
library contains 58 formal-search matches and 90 known-item supplements admitted
under common publication, version, cutoff, and evidence-role rules; 101 studies
receive the common classification. Because the candidate pool did not receive
complete full-text adjudication and the library is purposive, these counts
describe the review process rather than method frequency or deployment
prevalence.

The complete reported scheduling decision system is the conceptual unit of
analysis, and a consolidated substantive-work family is the bibliographic
counting unit. Materially different configurations within a family are
separated, and conclusions are stored as claim--evidence records. Direct
scheduling evidence requires an evaluated scheduling application in which a
learned artifact persists beyond the run that created it and affects a later
evaluated scheduling decision or its inference-time computation. Enabling and
adjacent evidence support mechanism interpretation and bounded research
directions, not scheduling-performance conclusions. EC.1--EC.6 document the
search, codebook, registries, evidence audit, and proposition tests.

\subsection{Organization of the Review}
\label{subsec:review-organization}

Sections~2--4 develop the taxonomy of inference architecture, learning and
adaptation, and decision assurance. Section~5 applies it to the evidence;
Section~6 develops the contingency framework and its proposed tests; and
Section~7 concludes.

\section{Inference Architecture: Learning Interfaces, Solution Formation, and Normal-Path Technical Authority}
\label{sec:decision_architectures}

This section classifies inference architecture through three analytically distinct questions: what a learned output does, whether and how learning forms schedule content inherited by the normal release path, and which components exercise binding control over releasable schedule content on that path. These questions must be recorded separately because similar learned predictions may determine a decoded schedule, provide an overwritable starting point, or only guide retained search.

The complete reported system is the conceptual unit of analysis, while each question uses its own coding unit, as Table~\ref{tab:architecture_records} shows. Conditional fallback and recovery are recorded under decision assurance in Section~4 rather than used to redefine normal-path technical authority. EC.2 provides the definitions and boundary rules, and EC.3 reports representative configurations and the 101-study common-field registry.

\subsection{Three Inference-Architecture Records}
\label{subsec:architecture_logic}

\begin{table}[!htbp]
\centering
\caption{Three records for inference architecture}
\label{tab:architecture_records}
\small
\begin{tabularx}{\textwidth}{P{2.5cm} P{4.2cm} Y P{3.7cm}}
\toprule
Record & Question & Coding unit & Primary values \\
\midrule
Learning interface
& What is the learned output used to do?
& Learned output--consumer link
& Surrogate substitution; candidate proposal; initialization; component
selection or configuration; search control; direct schedule formation \\
\addlinespace[2pt]
Solution formation
& Does learning form normal-path schedule content, and how?
& Normal release path; learned formation stage(s) when present
& Schedule-generation (SG) codes: SG0 no learned release-level formation; SG1 one-shot; SG2 constructive;
SG3 improvement; SG4 iterative global \\
\addlinespace[2pt]
Normal-path technical authority
& Who has binding control over schedule content?
& Complete normal release path
& Solver-led; shared technical authority; model-led \\
\bottomrule
\end{tabularx}
\end{table}

Table~\ref{tab:architecture_records} summarizes the three records, their coding units, and their primary values. The records are analytically distinct and can interact in a working system. For example, an SG1 stage may appear in either a system with shared technical authority or a model-led system. A learned warm start remains SG0 when the retained solver
can overwrite it, and a model can guide search without forming release-level
schedule content. Learning and adaptation are classified separately in
Section~3, and decision assurance is examined in Section~4. EC.2--EC.3 provide
the complete decision rules and mixed-pipeline cases.

\subsection{Learning Interfaces: What the Learned Output Does}
\label{subsec:learning_interfaces}

A learning interface is defined by the immediate downstream use of a learned output. Each output--consumer link receives a code, and materially different consumers are recorded separately. Similar predictions may replace a calculation, provide a candidate or starting state, guide retained search, or commit schedule content; the downstream consumer, rather than the output's form, determines the interface. These uses leave different downstream work and failure modes. Table~\ref{tab:interface_roles} presents the six interface families used in the review. EC.2 further divides search control into order-preserving and coverage-changing subcodes, yielding seven operational codes.

\begin{table}[!htbp]
\centering
\caption{Learning-interface families and representative scheduling examples}
\label{tab:interface_roles}
\small
\begin{tabularx}{\textwidth}{P{3.1cm} Y P{4.2cm}}
\toprule
Interface & Definition & Representative example \\
\midrule
Surrogate substitution
& A learned value replaces a numerical or logical calculation inside a retained procedure.
& Learned criterion values inside a decomposition heuristic
\citep{bouska2023deep}. \\
\addlinespace[2pt]
Candidate proposal
& Learning supplies an optional candidate object that a retained decision procedure may adopt, modify, or replace.
& Checked subproblem proposals inside a coordinated optimization method
\citep{liu2024integrating}. \\
\addlinespace[2pt]
Initialization
& Learning provides a starting solution or state that a retained procedure remains free to overwrite.
& Learned commitment schedules used to warm-start mixed-integer programming (MIP)
\citep{omalley2023hybrid}. \\
\addlinespace[2pt]
Component selection or configuration
& Learning selects a solver, rule, setting, operating mode, or handoff point before another procedure acts.
& Selection among constraint-programming (CP) solvers for flexible job shops
\citep{muller2022algorithm}. \\
\addlinespace[2pt]
Search control
& Learning orders or restricts alternatives inside retained search without directly forming release-level schedule content.
& Learned variable ordering inside complete constraint-programming search
\citep{sun2024enhancing}. \\
\addlinespace[2pt]
Direct schedule formation
& A learned output, possibly through a fixed decoder, commits a scheduling decision or relation inherited by the normal release path.
& One-shot priority sampling followed by fixed list scheduling
\citep{jeon2023neural}. \\
\bottomrule
\end{tabularx}
\end{table}

Interface boundaries follow the learned output's immediate downstream use. In \citet{bouska2023deep}, predicted values replace calculations inside a recursive scheduling heuristic. In \citet{liu2024integrating}, learning proposes subproblem solutions that retained optimization may coordinate or replace. The latter is therefore a candidate proposal because the downstream procedure determines whether and how the proposed content enters the final schedule.

Initialization and candidate proposal differ in their downstream roles. An initialization enters a retained procedure as its starting solution or state, whereas a proposal is an optional candidate evaluated or coordinated alongside other alternatives. Both remain nonbinding when the downstream procedure can replace them. Direct schedule formation applies when learned content is inherited by the normal release path and no downstream component can restore or replace it. A component that can only check or veto this content provides assurance. A retained optimizer may nevertheless determine additional schedule content, in which case the release path has shared technical authority.

Search control changes the order or coverage of alternatives while retained search constructs the schedule. Order-preserving guidance leaves the declared alternatives available; coverage-changing guidance may remove alternatives and thereby contribute binding control, even though learning does not directly form release-level content and the formation code remains SG0. Direct formation instead commits learned schedule content inherited by the normal release path. For example, learned fixing in rolling-horizon flexible job-shop scheduling commits selected machine relations before retained optimization determines the remaining assignments and timing \citep{li2025rolling}. The interface record therefore identifies the learned output's immediate use; solution formation and technical authority answer separate questions.

\subsection{Solution Formation: How a Schedule Is Produced}
\label{subsec:solution_formation}

The formation record states whether and how learning forms schedule content
inherited by the normal release path. A release-level object is a complete
candidate, a committed prefix, or a set of schedule relations inherited by the
completed schedule. SG0 records that no learned stage forms such an object.
SG1--SG4 describe four learned formation stages, as shown in
Figure~\ref{fig:generation_mechanisms}.

\begin{figure}[!htbp]
\centering
\begin{tikzpicture}[
    node distance=5mm and 6mm,
    every node/.style={font=\scriptsize},
    stagebox/.style={draw, rounded corners, align=center, minimum height=7mm, inner sep=2mm},
    stageflow/.style={->, thick}
]
  \node[stagebox] (i1) {instance};
  \node[stagebox, right=of i1] (m1) {simultaneous learned\\contribution};
  \node[stagebox, right=of m1] (d1) {complete or partial\\descriptor};
  \node[stagebox, right=of d1] (s1) {fixed decoding\\or preserving\\completion};
  \draw[stageflow] (i1)--(m1);
  \draw[stageflow] (m1)--(d1);
  \draw[stageflow] (d1)--(s1);
  \node[left=2mm of i1, anchor=east, font=\scriptsize\bfseries] {SG1 One-shot};

  \node[stagebox, below=8mm of i1] (i2) {instance and\\partial schedule};
  \node[stagebox, right=of i2] (m2) {conditional\\learned decision};
  \node[stagebox, right=of m2] (a2) {committed\\schedule content};
  \node[stagebox, right=of a2] (s2) {complete\\candidate};
  \draw[stageflow] (i2)--(m2);
  \draw[stageflow] (m2)--(a2);
  \draw[stageflow] (a2)--(s2);
  \draw[stageflow] (a2.south) |- ++(0,-4mm) -| (i2.south);
  \node[left=2mm of i2, anchor=east, font=\scriptsize\bfseries] {SG2 Constructive};

  \node[stagebox, below=8mm of i2] (i3) {complete candidate\\or incumbent};
  \node[stagebox, right=of i3] (m3) {learned move or\\search decision};
  \node[stagebox, right=of m3] (a3) {revised\\candidate};
  \node[stagebox, right=of a3] (s3) {returned\\candidate};
  \draw[stageflow] (i3)--(m3);
  \draw[stageflow] (m3)--(a3);
  \draw[stageflow] (a3)--(s3);
  \draw[stageflow] (a3.south) |- ++(0,-4mm) -| (i3.south);
  \node[left=2mm of i3, anchor=east, font=\scriptsize\bfseries] {SG3 Improvement};

  \node[stagebox, below=8mm of i3] (i4) {noise or relaxed\\global state};
  \node[stagebox, right=of i4] (m4) {learned global\\update};
  \node[stagebox, right=of m4] (a4) {updated\\global state};
  \node[stagebox, right=of a4] (s4) {decoded\\candidate};
  \draw[stageflow] (i4)--(m4);
  \draw[stageflow] (m4)--(a4);
  \draw[stageflow] (a4)--(s4);
  \draw[stageflow] (a4.south) |- ++(0,-4mm) -| (m4.south);
  \node[left=2mm of i4, anchor=east, align=right, font=\scriptsize\bfseries]
    {SG4 Iterative global};
\end{tikzpicture}
\caption{Four learned solution-formation stages. The codes describe how
learned schedule content changes during inference. Training and adaptation,
technical authority, and decision assurance are recorded separately. A system
may compose stages, such as SG1--SG3 or SG2--SG3.}
\label{fig:generation_mechanisms}
\end{figure}

Several widely used formation mechanisms originated in neural combinatorial optimization. Pointer Networks introduced variable-size sequential selection, and the Attention Model established an attention-based encoder--autoregressive-decoder template for constructive solution formation \citep{vinyals2015pointer,kool2019attention}. Policy Optimization with Multiple Optima (POMO) extended this constructive line through multiple-start policy optimization, while Matrix Encoding Network (MatNet) adapted attention to relation-matrix encoding and applied it directly to flexible flow-shop scheduling \citep{kwon2020pomo,kwon2021matnet}. GOAL and PolyNet further show how constructive architectures can be reused across tasks or solution strategies \citep{drakulic2025goal,hottung2025polynet}. Their scheduling applications provide direct evidence within the tested tasks; results from those tasks do not establish performance beyond the tested scheduling settings.

Graph neural networks (GNNs) and attention or Transformer encoders are
representation choices rather than solution-formation categories. Recent
Transformer-based systems report strong results against selected GNN-based
baselines \citep{choi2025multipolicy,xiao2026resched}. These comparisons also
change state representations, decoders, training procedures, or sampling
budgets. They therefore establish the performance of the tested configuration,
not a general ranking of Transformer and GNN encoders.

\paragraph{SG1: One-shot formation.}
An SG1 stage emits all learned schedule content at the same time. A fixed
decoder may turn that content into a feasible schedule, and a retained
procedure may complete or repair it while preserving binding learned content.
Neural DAG Scheduling, for example, operates on directed acyclic graphs (DAGs):
it samples all node priorities in one learned pass and applies list scheduling
\citep{jeon2023neural}.
\citet{kotary2022fast} predict job-shop start times in parallel; the predicted
times induce machine orders that are retained during normal-path feasibility
repair. Both contain an SG1 learned stage, although their technical-authority
configurations differ. A fixed decoder can enforce the constraints it encodes,
and an SG1 stage can be called again in a rolling-horizon system. Its evaluation
therefore includes the decoder, repair, and total release time.

\paragraph{SG2: Constructive formation.}
An SG2 stage forms schedule content through a sequence of commitments. Each
later learned decision depends on content committed earlier. Learned
dispatching policies for job-shop scheduling are standard examples
\citep{zhang2020dispatch,park2021learning}. MatNet also constructs flexible flow-shop schedules sequentially \citep{kwon2021matnet}. Constructive formation
describes inference, whereas reinforcement learning describes an update
mechanism. Constructive stages may also be learned from external targets or
other update mechanisms.

\paragraph{SG3: Improvement formation.}
An SG3 stage starts from a complete executable candidate and uses learned
decisions to revise it. The learned improvement heuristic of
\citet{zhang2024improvement} applies local moves to complete job-shop schedules,
while the Memory-enhanced Improvement Search framework (MIStar) combines a learned constructive initializer with a learned
improvement stage in its principal flexible-job-shop configuration
\citep{wang2025mistar}; EC.3 also records a variant with a nonlearned
initializer. Initialization, neighborhood coverage, and the full iteration
budget are part of the resulting system. A training procedure based
on self-improvement receives an SG3 inference code only when the reported
inference path revises a complete candidate.

\paragraph{SG4: Iterative global formation.}
An SG4 stage repeatedly updates a full-dimensional latent, noisy, or relaxed
representation. SG2 grows a partial schedule through successive commitments,
SG3 revises an executable incumbent, and SG4 updates a global representation
before decoding a candidate. Diffusion solvers such as Graph-Based Diffusion Solvers for Combinatorial Optimization (DIFUSCO) provide adjacent
combinatorial-optimization evidence for this mechanism \citep{sun2023difusco}.
Direct scheduling evidence is emerging, with recent work applying diffusion to
flow-shop and job-shop settings \citep{li2026goaldiffusion}. We therefore treat
SG4 as an emerging scheduling category.

Two boundary cases clarify the coding rule. Repeated model calls receive an SG2
code only when later learned decisions condition on earlier committed schedule
content. One model call receives an SG1 code only when its output forms schedule
content inherited by the normal path. DeepACO produces learned heuristic
measures that guide ant-colony optimization (ACO); because retained search forms the
schedule, it is coded as search control with SG0 formation
\citep{ye2023deepaco}. Systems may also compose stages. For example, a one-shot
or constructive stage may produce a schedule that a learned improvement stage
then revises, giving SG1--SG3 or SG2--SG3.

The formation record identifies the temporal structure of learned
schedule-forming inference and the model calls, samples, and iterations that
belong in full-pipeline runtime. Training, assurance, and technical authority
are recorded separately.

\subsection{Normal-Path Technical Authority: Who Controls Schedule Content?}
\label{subsec:decision_authority}

Technical authority is coded for a complete normal release path. A learned
output has a \emph{binding schedule effect} when it fixes or removes schedule
choices and no later component on that path retains the structural right and
decision space needed to restore or replace them. Organizational rights over
specification, release, intervention, and system lifecycle are analyzed
separately in Section~6.

The classification follows two questions. First, does any learned output have
a binding schedule effect on the normal release path? If not, the path is
\emph{solver-led}. If so, does a retained nonlearned schedule-forming procedure
independently determine additional binding schedule content on the same path?
If yes, the path has \emph{shared technical authority}; otherwise, it is
\emph{model-led}. EC.2 gives the complete test for multiple reported variants.

In a \emph{solver-led} path, learning provides advice, routing, or an initial
state, while a retained procedure remains free to form the schedule. Learned
variable ordering in complete constraint-programming search is solver-led when
the search preserves the required branches \citep{sun2024enhancing}. The hybrid
configuration of \citet{omalley2023hybrid} is also solver-led because the
mixed-integer solver can overwrite the learned warm starts and retains the
original decision space.

In a \emph{model-led} path, learned output fixes consequential schedule
content, while nonlearned components apply fixed decoding, constraint
propagation, checking, or veto without independently selecting alternative
schedule content. Neural DAG Scheduling is model-led because learned priorities
determine the list-scheduling result \citep{jeon2023neural}.

Under \emph{shared technical authority}, learned and retained schedule-forming
procedures both make binding contributions on the normal path. In
\citet{kotary2022fast}, learned predictions determine machine orders that are
preserved during feasibility repair, while retained optimization determines
the final timing. Learning-guided rolling-horizon optimization also has shared
technical authority when learned fixing commits some machine relations and the
retained optimizer chooses the remaining schedule content
\citep{li2025rolling}.

Three boundaries complete the classification. The presence of a solver does
not make a path solver-led after learning has fixed schedule choices that the
solver cannot restore. A checker or veto supplies decision assurance when it
can only permit or block release. A conditional fallback is recorded under
recovery assurance and does not alter normal-path technical authority.
Solver-led paths, paths with shared technical authority, and model-led paths
are alternative system configurations rather than stages of technical maturity.

\subsection{System Assembly and Comparative Implications}
\label{subsec:system_assembly}

The three records must be read together. Table~\ref{tab:architecture_examples} shows that the same formation stage can appear under different technical-authority configurations and that SG0 does not by itself imply a solver-led path. The table provides representative cases; EC.3 reports more detailed configuration records and the 101-study common-field registry.

\begin{table}[!htbp]
\centering
\caption{Representative combinations of interface, formation, and normal-path technical authority}
\label{tab:architecture_examples}
\small
\begin{tabularx}{\textwidth}{Y P{3.1cm} P{2.6cm} P{2.8cm}}
\toprule
Representative configuration & Interface & Formation & Technical authority \\
\midrule
Learned ordering within complete CP search \citep{sun2024enhancing}
& Search control (order-preserving) & SG0 & Solver-led \\
Learned schedules used as MIP warm starts \citep{omalley2023hybrid}
& Initialization & SG0 & Solver-led \\
Learned action-set restriction followed by retained scenario evaluation
\citep{omalley2023hybrid}
& Search control (coverage-changing) & SG0 & Shared technical authority \\
Parallel prediction followed by normal-path linear programming (LP) completion
\citep{kotary2022fast}
& Direct schedule formation & SG1 & Shared technical authority \\
One-shot priorities followed by fixed list scheduling \citep{jeon2023neural}
& Direct schedule formation & SG1 & Model-led \\
Constructive learned dispatching
\citep{zhang2020dispatch,park2021learning}
& Direct schedule formation & SG2 & Model-led \\
Learned improvement of a dispatching-rule schedule
\citep{zhang2024improvement}
& Direct schedule formation & SG3 & Shared technical authority \\
Neural heuristic measures within ACO \citep{ye2023deepaco}
& Search control (order-preserving) & SG0 & Solver-led \\
\bottomrule
\end{tabularx}
\end{table}

These combinations show why the three records should not be collapsed. Inference architecture follows the downstream decision process rather than a model-family name. The same SG1 stage can appear in a shared-authority or model-led path, while SG0 can appear in a solver-led or shared-authority path depending on whether search guidance preserves or removes alternatives. Technical authority therefore depends on whether learned effects are binding and whether a retained procedure independently contributes additional binding schedule content.

Assurance requirements are determined separately by the scope of binding learned control, the consequences of error, uncertainty about applicability, and the availability of recovery before the decision deadline. Sections~3 and~4 next examine the information used to create or update learned capability and the controls that support release.

\section{Learning and Adaptation Signals: Information Infrastructure and Operational Consequences}
\label{sec:training_signals}

Section~\ref{sec:decision_architectures} classifies the deployed inference
path. This section separately classifies the task information used to create
or update learned capability. The distinction matters because the same
constructive decoder may learn from solver targets, selected model candidates,
or sampled returns, while the same signal may support solver guidance or direct
schedule formation.

We code each separable update term or phase by the information that directly
defines it. \emph{Reference-target learning} fits a target supplied
independently of the current generator or selected from its candidates.
\emph{Direct-evaluation learning} lets a differentiable problem loss or sampled
return or fitness define the update. These branches yield the four mechanisms
in Table~\ref{tab:training_signal_comparison}. Studies may combine them
concurrently or sequentially; EC.2 provides the detailed coding rules.

\begin{table}[t]
\centering
\caption{Learning and adaptation signals impose different information and
validation requirements}
\label{tab:training_signal_comparison}
\begingroup
\setlength{\tabcolsep}{3pt}
\renewcommand{\arraystretch}{1.10}
\scriptsize
\begin{tabularx}{\textwidth}{P{2.55cm} P{4.20cm} P{3.25cm} Y}
\toprule
Update mechanism
& Operational coding rule
& Required information infrastructure
& Principal validity risk \\
\midrule

Independent-reference targets
& Fit a target that is treated as fixed during the update and does not depend
on candidates produced by the current learner; examples include solver,
heuristic, human, logged-action, and fixed-teacher targets.
& A renewable target-generation and audit pipeline
& Teacher bias, target cost, limited state coverage, and ambiguity among
equivalent schedules \\

Model-generated selected targets
& Generate candidates, evaluate and select or accept them, and then fit the
selected candidates as targets held fixed during the update.
& Diverse candidate generation, a faithful evaluator, and an acceptance rule
& Self-confirmation, inadequate diversity, objective mismatch, and hidden
rollout cost \\

Differentiable problem-loss updates
& Propagate a gradient from an operational objective, constraint, or explicit
surrogate through a differentiable decision interface to the updated variables.
& A credible differentiable interface and a valid discretization or decoder
& Proxy or relaxation error, invalid discrete outputs, and unreported
adaptation cost \\

Return- or fitness-based updates
& Use sampled or logged rewards, returns, advantages, Bellman targets, or
fitness to update a policy, value function, or parameterized distribution.
& An interactive environment, an objective evaluator, or logged experience
with adequate support
& Reward misspecification, simulator error, weak data coverage, and expensive
interaction \\

\bottomrule
\end{tabularx}

\vspace{2pt}
\parbox{\textwidth}{\scriptsize
\emph{Notes.}
Each row codes the task information that directly defines an update. Pretext
relations defined independently of the current learner are auxiliary reference
signals. Generic regularizers are recorded separately from task-informing
updates. Data acquisition, update timing, update object, and signal composition
are recorded as modifiers. Online reinforcement learning (RL) denotes
interactive experience acquisition; live production learning is recorded
separately.}
\endgroup
\end{table}

\subsection{Reference-Target Learning}
\label{subsec:reference_target_learning}

\emph{Independent-reference target learning} fits actions, schedules, solver
quantities, or relations supplied independently of the current model's
candidates. Scheduling studies have learned priority rules from optimization
solutions, queried an oracle on learner-visited states, and predicted quantities
used by a retained algorithm
\citep{ingimundardottir2011supervised,ingimundardottir2018discovering,
parmentier2023structured}. Adaptive oracle queries remain in this category: the
learner chooses the state, but the oracle supplies the target. Target provenance
should be stated precisely. For example, \citet{li2024imitation} use an improved
Nawaz--Enscore--Ham heuristic as the teacher for permutation flow-shop
scheduling.

The value of this mechanism depends on the target pipeline. Exact targets can
be expensive as instances and constraints grow; heuristic or historical targets
can reproduce existing weaknesses; and equivalent schedules can make
sequence-level losses misleading. Evidence should therefore report the target
producer, solve status or teacher quality, generation cost, coverage, and
treatment of alternative optima. Supervised learning can be comparatively straightforward to operationalize when a credible and renewable target-generation process already exists.

\emph{Model-generated selected-target learning} generates candidates,
evaluates them, selects or accepts some of them, and fits the selected
candidates as fixed targets. The Self-Labeling Improvement Method (SLIM) uses the best evaluated model sample as a
pseudo-label \citep{corsini2024selflabeling}. Later methods improve candidate
diversity through multiple policies, larger sample sets, or scheduling
symmetries
\citep{pirnay2024selfimprovement,choi2025multipolicy,
luttmann2026multiaction}. Methods described as \emph{bootstrapping} or
\emph{self-improvement} enter this category when the current model generates
the candidates used as later targets. We use \emph{self-labeling} for this
mechanism; the broader term \emph{self-supervised learning} also includes
representation-level pretext tasks.

Self-labeling removes the external teacher but still requires candidate
generation and a faithful evaluator. Evidence for self-labeling is most credible when rollouts are sufficiently inexpensive and diverse, candidate quality can be evaluated accurately, and the acceptance rule limits reinforcement of weak candidates.
Studies should report candidate and objective calls, acceptance rates,
diversity, and total training time. Independent targets concentrate risk in
the teacher; model-generated targets concentrate it in evaluation, diversity,
and the acceptance rule.

\subsection{Direct-Evaluation Learning}
\label{subsec:direct_evaluation_learning}

Self-labeling uses evaluation to select a fixed target. Direct-evaluation
methods instead let the evaluation itself define the update.

\emph{Differentiable problem-loss updates} propagate a gradient from an
objective, constraint, or explicit surrogate through a differentiable decision
interface. This mechanism is distinct from ranking candidates by makespan,
which creates selected targets, and from using makespan as an RL reward, which
creates return-based updates. Mixed losses may combine independent targets
with a differentiable constraint term; inference-time feasibility still
depends on the decoder, completion, and checking path
\citep{kotary2022fast}. Test-time adaptation through a feasibility surrogate is
an ordered problem-loss phase following supervised initialization
\citep{liu2025integrated}.

Reusable scheduling evidence for this mechanism remains concentrated in such
mixed systems. Differentiable Combinatorial Scheduling at Scale instead
optimizes relaxed decision variables separately for each instance
\citep{liu2024differentiable}. It demonstrates problem-loss updating, but does
not create an artifact reused on later instances. This boundary separates a
differentiable per-instance optimizer from a reusable learned scheduler.

\emph{Return- or fitness-based learning} updates a policy, value function, or
parameter distribution from returns, temporal-difference targets, advantages,
or black-box fitness. Scheduling applications precede modern graph policies
\citep{zhang1995rljobshop,riedmiller1999neural,aydin2000dynamic,
gabel2012distributed}. Recent studies use graph-based constructive policies
\citep{zhang2020dispatch,park2021learning,song2023flexible} and
attention- or Transformer-based designs
\citep{kwon2021matnet,wang2024dualattention,xiao2026resched}.
The Attention Model established an influential REINFORCE-based training design
for constructive neural optimization \citep{kool2019attention}. POMO generates
trajectories from multiple starting points for each instance and uses their
mean return as a shared REINFORCE baseline, reducing gradient variance
\citep{kwon2020pomo}. MatNet applies this NCO line directly to flexible flow-shop scheduling \citep{kwon2021matnet}. Recent equipment-specific work also uses REINFORCE with a shared baseline for noncyclic cluster-tool scheduling \citep{kim2025noncyclic}.
Black-box fitness can also update a parameter distribution; \citet{kaleta2026constructive}, for example, use the covariance matrix adaptation evolution strategy to update neural-network weights from makespan-based fitness. These examples reinforce
the distinction between constructive inference and return- or fitness-based
learning.

Online and offline RL are acquisition modifiers. Online RL obtains experience
through interaction; offline RL learns from fixed transitions and rewards, as
in offline learning to dispatch \citep{vanremmerden2025offline}. Logged actions
used as targets instead constitute imitation. Self-critical training and POMO
remain return-based because their comparison rollouts provide baselines, not
targets.

Terminal rewards preserve direct contact with the reported objective but delay
credit. Shaped rewards provide denser information and require evidence that
they remain aligned with that objective
\citep{ng1999reward,zhang2020dispatch}. The scheduling literature contains both proximal policy optimization (PPO) and REINFORCE systems; optimizer choice alone does not determine inference
architecture or operational suitability.

Dynamic and stochastic scheduling are also distinct. Dynamic scheduling
reveals information or revises decisions during operation, as with new-job
insertion or changing workflow arrivals
\citep{luo2020dynamic,yang2025goodrl}. Stochastic scheduling represents
uncertain quantities or events and may produce either an ex ante schedule or an
adaptive policy
\citep{smit2025stochastic,zhang2026uncertainty}. Return-based learning is most
credible when the state captures relevant information, rewards match the
operating objective, and the simulator or logged data cover consequential
state--action choices.

\subsection{Signal Composition and Operational Implications}
\label{subsec:signal_composition}

Systems may combine update mechanisms. We code simultaneous loss terms as
\emph{concurrent} composition and successive phases as \emph{ordered}
composition. Solver labels plus a differentiable constraint loss form a
concurrent independent-reference/problem-loss regime. Independent-reference pretraining
followed by instance-level surrogate adaptation forms an ordered regime
\citep{liu2025integrated}. Heterogeneous Earliest Finish Time (HEFT) imitation followed by simulator-interactive PPO
and simulated operating-time PPO adaptation forms an ordered
independent-reference/return-based regime with a later adaptation phase
\citep{yang2025goodrl}. Component ablations should separate the value of signal
composition from the value of additional labels, objective calls, solver
information, or computation.

Four modifiers are coded separately: how information is acquired, when updates occur, what variables are updated, and how mechanisms are combined. Fixed logs may support either imitation or offline RL; test-time describes when an update occurs rather than what defines it; and simulator interaction remains distinct from live production learning.

Each mechanism relies on a different information resource. Independent targets
require a teacher; model-generated targets require diverse candidates and a
faithful evaluator; differentiable losses require a valid gradient interface;
and return- or fitness-based learning requires a simulator, supported logs, or
an objective evaluator. Label-free methods therefore replace target-generation
cost with other infrastructure. Learning information explains how capability is
acquired; release still depends on feasibility, performance, applicability, and
recovery. Section~\ref{sec:decision_assurance} examines those tasks,
Section~\ref{sec:critical_evidence} audits their cost and validity, and
Section~\ref{sec:om_contingency} treats the required infrastructure as an
organizational contingency.

\section{Decision Assurance: Feasibility, Performance, Applicability, and Recovery}
\label{sec:decision_assurance}

Sections~\ref{sec:decision_architectures} and
\ref{sec:training_signals} classify how reported systems form schedules and
acquire learned capability. Across the audited analytical library, we did not identify a common vocabulary for the controls between model output and release. We synthesize those controls
through a 3+1 framework: \emph{feasibility}, \emph{performance}, and
\emph{applicability} are three claims supporting release, while
\emph{recovery} is the conditional response when the normal path fails a check
or misses its deadline. The framework provides common records for comparing
heterogeneous studies.

Each record links a required outcome to the reported control, supporting
evidence, scope, and failure response. Scope includes the represented
constraints, problem and information regimes, operating conditions,
tolerances, and period of validity. For candidates produced through recovery,
we separately record the reported support for the three release claims.
Normal-path technical authority remains separate from organizational rights
over release and intervention. Table~\ref{tab:assurance_mechanisms} gives
representative source records; EC.2 provides the coding rules, and EC.3--EC.5
report the fuller evidence base.

\begin{table}[t]
\centering
\caption{Decision-assurance records used in the review}
\label{tab:assurance_mechanisms}
\begingroup
\setlength{\tabcolsep}{3pt}
\renewcommand{\arraystretch}{1.10}
\scriptsize
\begin{tabularx}{\textwidth}{P{1.85cm} P{3.35cm} P{3.55cm} Y}
\toprule
Record
& Review question and coding unit
& Control types coded
& Representative evidence and boundary \\
\midrule

Feasibility
& Does the final candidate satisfy the stated constraints? \newline
Candidate and complete formation--checking path
& Constraint-preserving formation; retained LP or solver completion;
independent final-candidate checking (unresolved in the representative records)
& Fixed decoding and eligible-action construction preserve represented
constraints \citep{jeon2023neural,zhang2020dispatch}; retained optimization
completes restricted formulations \citep{kotary2022fast,li2025rolling}.
Claims remain constraint-specific. \\

Performance
& What quality or risk statement is supported? \newline
Candidate or policy under a stated formulation or event model
& Direct objective or key performance indicator (KPI) evaluation; valid bound or proof status;
statistical evaluation
& Problem-specific bounds and retained-solver proof status provide
deterministic evidence
\citep{muller2022algorithm,omalley2023hybrid,sun2024enhancing}; scenario
evaluation supports expected or risk performance
\citep{omalley2023hybrid,smit2025stochastic}. \\

Applicability
& Does the evidence cover the current instance and operating conditions?
\newline Evidence--case link
& Bounded validation; staged introduction; monitoring; revalidation when
reported
& \citet{stockermann2025deploymentreview} identify validation, transfer, and
monitoring concerns; \citet{konstantelos2022fabwide} report staged validation
and monitoring; \citet{kim2024photolithography} provide bounded field-use
evidence. \\

Recovery
& What follows a failed check or missed deadline? \newline
Conditional failure episode and alternative path
& Conditional greedy or solver backup; planner intervention
& The literature reports an untriggered greedy branch
\citep{kotary2022fast}, an exercised internal solver backup
\citep{liu2024integrating}, and planner intervention
\citep{liang2022lenovo}. These records cover different system levels. \\

\bottomrule
\end{tabularx}

\vspace{2pt}
\parbox{\textwidth}{\scriptsize
\emph{Notes.}
Feasibility, performance, and applicability are nonexclusive release claims;
recovery is a conditional path. The examples show how the records are applied
and do not estimate the prevalence of assurance controls.}
\endgroup
\end{table}

\subsection{Feasibility: What the Reported Controls Establish}
\label{subsec:establish_restore_feasibility}

The reviewed studies establish feasibility in two main ways. Some embed the
constraint logic in schedule formation. One-shot priorities followed by fixed
list scheduling and constructive selection among eligible operations satisfy
the precedence and resource constraints represented by their respective
decision rules \citep{jeon2023neural,zhang2020dispatch}. Feasibility in these
systems belongs to the complete learned-output--decoder path, not to the
one-shot or constructive label.

Other studies retain optimization after prediction. In
\citet{li2025rolling}, a learned classifier selects overlapping operations
whose previous-horizon machine assignments are then fixed; the CP-SAT solver, which combines constraint programming and Boolean satisfiability, determines
the remaining assignments and timing. In \citet{kotary2022fast}, predicted
start times induce a machine order and a normal-path linear program supplies
feasible timing when that order is consistent. These are
formulation-specific completion results, not general feasibility or optimality
certificates.

Constraint coverage is therefore the decisive boundary. Scheduling settings
may add routing flexibility, assembly relations, sequence-dependent setups,
and online arrivals to basic precedence and capacity constraints
\citep{reijnen2026benchmark}. Reported masks, decoders, and completion models
cover only what they encode. An independently implemented final-candidate checker would also belong to feasibility because it verifies the completed schedule against the represented constraints and tolerances after formation. Direct evidence for independently implemented final-candidate checking is absent from these representative records.
\citet{konstantelos2022fabwide} report trial periods and automatic validation,
but the source does not establish a candidate-level independent checker; we
therefore leave that control unresolved.

\subsection{Performance: Bounds and Statistical Evidence}
\label{subsec:quality_performance}

The literature supports three different performance statements. Direct
objective or KPI evaluation can rank candidates or apply a stated selection
rule, as in the parallel comparison of \citet{kim2024photolithography}. A
deterministic certificate instead requires a feasible candidate and a valid
bound or proof status for the same formulation. Under uncertainty, performance
is estimated for a specified event model and risk measure.

Retained solvers provide the clearest deterministic evidence. In
\citet{muller2022algorithm}, learning selects a CP solver. The study evaluates
its returned schedule against a problem-specific lower bound and separately
records whether the selected solver proves optimality. In
\citet{sun2024enhancing}, learned ordering changes search effort, but the
retained solver establishes optimality only upon completed proof. Time-limited
incumbents and comparator differences support empirical performance statements
rather than certificates.

The stochastic studies expose a second boundary. \citet{omalley2023hybrid}
separate the MIP gap for a finite-scenario unit-commitment formulation from
Monte Carlo estimates of operating cost. \citet{smit2025stochastic} evaluate
expected makespan and value at risk under specified processing-time
distributions and a fixed scenario budget. These results support the reported
formulations, distributions, and evaluation budgets. They also show why
performance evidence and the applicability of that evidence require separate
records.

\subsection{Applicability: The Scope of Evidence}
\label{subsec:quality_applicability}

Applicability concerns whether the evidence supporting release covers the
current decision. Relevant boundaries include problem structure and scale,
decision-time information and objectives, event processes, tested shifts, and
the period over which validation remains informative. The deployment review of
\citet{stockermann2025deploymentreview} identifies validation,
simulator-to-field transfer, and monitoring as recurring concerns, providing a
domain-level basis for this record.

The field evidence remains setting-specific. In the solver-led deployment of
\citet{konstantelos2022fabwide}, localized trials, adjustable guidance
strength, automatic validation, and operator monitoring precede broader
control. \citet{kim2024photolithography} document one year of
actual-fabrication use, with published comparisons for four consecutive
months. These studies establish operational possibility and bounded validation
in their reported fabs. Across the coded studies, direct evidence on formal
boundary detection, recalibration, and post-change revalidation is scarcer than
evidence on feasibility or performance.

\subsection{Recovery: Evidence on Conditional Response}
\label{subsec:release_recovery}

Recovery begins when the normal path fails a required check or cannot produce
a releasable candidate by the deadline. It may occur inside an algorithm,
through an alternative scheduling path, or through operational intervention.
For each resulting candidate, the review records any reported feasibility,
performance, and applicability evidence separately.

The reviewed evidence spans these levels. \citet{kotary2022fast} define a
greedy procedure for an inconsistent learned machine order, but the branch was
not needed in their tests. \citet{liu2024integrating} exercise branch-and-cut
when a learned part-subproblem proposal fails its acceptance condition,
providing evidence for an internal algorithmic backup rather than a complete
release-level recovery process. \citet{liang2022lenovo} report planner review,
manual modification, and replanning in operational use, but do not quantify
recovery frequency, capacity, or deadline success.

Other controls mark useful boundaries. The parallel learned and rule-based
candidates in \citet{kim2024photolithography} belong to normal-path selection,
while the staged rollout in \citet{konstantelos2022fabwide} supports
applicability and release governance. Neither record documents
failure-triggered recovery. Across the coded records in EC.3--EC.5, fallback
mechanisms and human intervention are more visible than evidence on complete,
timely recovery. Section~\ref{sec:critical_evidence} carries this gap into the
evidence audit, and Section~\ref{sec:om_contingency} examines its implications
for technical authority and organizational decision rights.

\section{Applying the System Taxonomy: What the Audited Evidence Supports}
\label{sec:critical_evidence}

Sections~\ref{sec:decision_architectures}--\ref{sec:decision_assurance}
provide records for inference architecture, learning and adaptation, and
decision assurance. This section applies those records to four evidence
questions: system performance and cost, transfer of learned capability,
applicability under operating variation, and operational use. The unit of
assessment is a claim about a reported configuration paired with the evidence
relevant to that claim. EC.3--EC.5 provide the configuration records and the
claim--evidence records for technical and application claims.

The taxonomy is not used to rank solver-led configurations, configurations
with shared technical authority, and model-led configurations. It determines
whether reported results concern comparable decision objects, cost boundaries,
guarantees, operating changes, and release conditions. The architecture records
identify the work retained after prediction and the decision space controlled
by each component. The learning record distinguishes transfer of retained
capability from retraining or architectural reuse, while the assurance and
maturity records separate operating use from automated release.

Applying this frame changes the interpretation of several familiar claims.
Similar model results may leave different downstream work, and a solver's
guarantee covers only the formulation and decision space it is allowed to
control. Learning may reduce recurring computation while moving cost to
capability acquisition, assurance, adaptation, or maintenance. Table~\ref{tab:claim_evidence_summary}
shows how these distinctions structure the audit and which comparisons remain
unresolved.

\begin{table}[!htbp]
\centering
\caption{How the system taxonomy structures the evidence audit}
\label{tab:claim_evidence_summary}
\begingroup
\setlength{\tabcolsep}{3pt}
\renewcommand{\arraystretch}{1.08}
\scriptsize
\begin{tabularx}{\textwidth}{P{2.45cm} P{4.05cm} Y}
\toprule
Evidence question
& Taxonomy records applied
& System-level conclusion and unresolved comparison \\
\midrule

Performance and total cost
& Learning interface, solution formation, normal-path technical authority,
update mechanisms, and assurance
& Configurations may differ in downstream work, recurring cost, and the
decision space covered by a guarantee. Complete-path comparisons under a common
release contract and lifecycle boundary remain unresolved. \\

Transfer of learned capability
& Retained learned artifact, train--test change, update phase, adaptation work,
and applicability
& Transfer requires retained capability, a specified change, and an adaptation
budget. Evidence beyond within-family size changes remains bounded. \\

Operating variation
& Information and objective regime, update timing, and applicability
& Performance within a declared event, uncertainty, or preference regime does
not establish adaptation to a different regime or scheduling structure. \\

Operational use
& Operational maturity, technical authority, assurance and recovery, and
organizational rights
& Contact with operating decisions does not by itself establish binding learned
control, automated release, or recovery capacity. Their joint effects have not
been compared on a common basis. \\

\bottomrule
\end{tabularx}

\vspace{2pt}
\parbox{\textwidth}{\scriptsize
\emph{Notes.} The records determine the scope of each claim. The conclusions
refer to the audited evidence and do not estimate the frequency of methods or
operating configurations.}
\endgroup
\end{table}

\subsection{System-Equivalent Performance and Total Decision Cost}
\label{subsec:benchmarking_quality}

A system-equivalent comparison holds the scheduling formulation and available
information, release requirement, computing environment, and decision budget
sufficiently constant. The architecture records determine what enters that
budget. SG0 search guidance leaves search downstream, while direct formation
may still require decoding, repair, checking, or recovery. Model inference time
is therefore not a common measure of complete-system response time. Shared
benchmarks align instances
\citep{taillard1993benchmarks,brandimarte1993routing}, and unified environments
can reduce formulation and implementation differences
\citep{reijnen2026benchmark}.

The reviewed comparisons are informative but support different claims. Learning to Dispatch compares its learned rule with priority rules and uses Google
OR-Tools solutions obtained under a 3,600-second limit, or best-known values, as quality references \citep{zhang2020dispatch}. \citet{zhao2024heterogeneous} compare their heterogeneous-graph policy with dispatching rules and a CPLEX run capped at 1,800 seconds. When CPLEX does not prove optimality within that limit, its incumbent is nevertheless used as the reference value. The study also reports comparisons with an adapted double deep Q-network and a genetic algorithm on a separate extra-large test set.

These experiments establish performance relative to their reported comparison sets, not state-of-the-art performance under a common release contract.
Comparator coverage and decision budgets differ, and the CPLEX reference is not always an optimality certificate. Imitation learning against an improved
Nawaz--Enscore--Ham teacher provides a more direct size-dependent quality--computation-time comparison, but that result remains bounded by the
teacher, tested sizes, and reported computation \citep{li2024imitation}.
Taken together, the evidence supports configuration-specific trade-offs rather than a general ranking of learning, exact optimization, and metaheuristics.

Technical authority also localizes a reported guarantee. With order-preserving
guidance, complete constraint-programming search can prove optimality for the
declared formulation when the proof terminates \citep{sun2024enhancing}. When
learned output fixes machine relations that subsequent optimization cannot
reopen, retained optimization establishes feasibility or performance only
within the resulting decision space \citep{kotary2022fast,li2025rolling}. The
presence of a solver therefore does not make guarantees equivalent across
systems.

Learning also reallocates rather than removes decision work. Recurring cost
includes formation, sampling, retained search, assurance, and conditional
response. Lifecycle cost additionally includes targets or simulation,
training, integration, updating, and maintenance. Self-labeling moves work to
candidate generation and evaluation during capability acquisition
\citep{corsini2024selflabeling}; iterative improvement and rolling-horizon
configurations add repeated evaluation or solver calls during use
\citep{wang2025mistar,li2025rolling}. The audited evidence supports selected
quality--computation-time trade-offs, but does not provide a common basis for
complete-pipeline lifecycle economics.

\subsection{Transfer of Learned Capability}
\label{subsec:scalability_generalization}

Transfer is assessed through the learned artifact retained from development, a
specified change, the adaptation performed and its budget, and the evidence
established afterward. Computational reach is a weaker claim: it shows that a
method executes at a new scale. Performance at scale additionally requires
reported feasibility, quality, and decision time. Reusing an architecture after
separate training establishes architectural reuse rather than transfer of
learned capability.

Direct scheduling evidence supports bounded within-family size transfer.
The cited studies retain trained artifacts and evaluate them on larger instances
within the same scheduling family
\citep{zhang2020dispatch,park2021learning,li2024imitation,xiao2026resched}.
Their claims are defined by the reported size ranges, comparators, release
paths, and adaptation budgets. Tests under changed processing-time
distributions provide further evidence within stated problem and information
regimes \citep{park2021learning,smit2025stochastic}.

Related vehicle-routing studies examine multi-task transfer through
mixture-of-experts and knowledge distillation
\citep{zhou2024mvmoe,zheng2025mtlkd}. Scale-diverse training in the related
container-retrieval problem provides an adjacent cross-size design \citep{shin2026containers}. These studies motivate analogous cross-task and
cross-size tests in scheduling but do not establish transfer of learned
scheduling capability.
Changes in formulation or task require a different interpretation. RESCHED
uses a common architecture but trains separate models for the flexible job-shop scheduling problem (FJSP), the job-shop scheduling problem (JSSP), and the flexible flow-shop scheduling problem (FFSP)
\citep{xiao2026resched}. GOAL instead reports adaptation-based cross-problem
transfer by fine-tuning a multitask backbone trained partly on JSSP for
open-shop scheduling \citep{drakulic2025goal}. PolyNet generates diverse
strategies within FFSP rather than transferring a learned artifact across
scheduling tasks \citep{hottung2025polynet}. These distinctions make transfer
a system-level claim about retained capability, adaptation work, and renewed
assurance. The direct evidence supports stated size ranges, selected parameter
changes, and the reported cross-problem adaptation; it does not establish
general one-policy transfer across scheduling structures or tasks.

\subsection{Applicability under Operating Change, Uncertainty, and Preferences}
\label{subsec:changing_conditions}

Dynamic scheduling evidence concerns a declared information-arrival and
replanning process. Studies with new-job arrivals, workflow arrivals, and
rolling-horizon events report quality and decision-time outcomes under their
modeled processes \citep{luo2020dynamic,yang2025goodrl,li2025rolling}. The
update record matters: repeated scheduling with a fixed policy demonstrates
response within the regime, whereas retraining or test-time updating adapts the
learned capability. Dynamic operation alone does not establish adaptation to a
new event process.

Stochastic evidence is bounded by the uncertainty model, the information
revealed before each decision, and the reported performance measure.
Scenario-based flexible-job-shop scheduling has been evaluated for expected
makespan and value at risk under specified processing-time distributions
\citep{smit2025stochastic}. Work on structural uncertainty optimizes expected
performance and reports conditional value at risk (CVaR) as an evaluation measure
\citep{zhang2026uncertainty}. These studies establish performance for their
stated stochastic models rather than a general response to unmodeled change.

Preference variation changes another part of the applicability boundary.
\citet{luo2021multiobjective} switch among scalar objectives in a dynamic
setting; \citet{ding2024multipolicy} use multiple policies to cover weight
settings and generate Pareto trade-offs; and
\citet{smit2026multiobjective} test a preference-conditioned policy on
untrained weights within a stated domain. These records support trade-off
generation and bounded preference generalization. Changes to constraints,
routing, resource relations, or the formulation itself require separate
adaptation and renewed assurance. Across the three streams, the evidence
supports applicability within declared operating regimes, not general
robustness across regimes.

\subsection{Operational Maturity, Authority, and Release}
\label{subsec:industrial_evidence}

Operational maturity records the closest documented contact with operating
decisions and remains separate from technical authority. Level A requires
documented routine or sustained use in an ongoing process. Level B1 covers a
bounded physical experiment, pilot, shadow operation, or controlled testbed
without documented routine use. EC.5 distinguishes these records from
retrospective operational-data evaluation, industry-calibrated simulation, and
synthetic evaluation.

The source-adjudicated register contains four learning-augmented Level A
episodes, one solver-led Level A industrial reference case, and five
learning-augmented Level B1 experiments or testbeds. The Level A
learning-augmented cases cover assembly planning, photolithography,
scheduling for naphtha cracking, and fab scheduling
\citep{liang2022lenovo,kim2024photolithography,hong2024naphtha,sood2024fab}.
Lenovo reports rollout across four plants and 43 lines together with planner
review and modification. The photolithography study documents one year of
actual-fabrication use under its reported release arrangement. The solver-led
Seagate case reports continuous operation from June 2022, preceded by staged
validation and accompanied by operator monitoring
\citep{konstantelos2022fabwide}.

These cases show why maturity, authority, and release must be coded separately. Level A establishes operating use, but the binding contributions of learned and retained components still determine technical authority. Automated release depends on the reported gate and organizational release rights; fallback and intervention require separate recovery evidence. The five Level B1 records likewise establish bounded end-to-end execution---one physical manufacturing experiment and four computing testbeds---without documenting routine use \citep{zhang2026digitaltwin,narayanan2020gavel,mao2019decima, jajoo2022slearn,delimitrou2013paragon}. Beyond these Level A and B1 records, learning-augmented scheduling has also been evaluated for noncyclic and concurrent cluster tools and for robotic-cell and robotic-flow-shop settings
\citep{kim2025noncyclic,kim2025concurrent,kim2024roboticcell, kim2022dualgripper,lee2022roboticflowshop,kim2023lookahead}.
These studies broaden the equipment-specific scheduling evidence, including look-ahead search and processing-time variation, but report computational evaluation rather than documented routine use or bounded physical execution. EC.1 records their individual evidence boundaries.

The application evidence establishes operational possibility for different
configurations within the duration, release, intervention, maintenance, and
rollback record reported for each case. This reading accords with deployment
reviews that identify validation, field transfer, and monitoring as continuing
concerns \citep{stockermann2025deploymentreview}.

Applying the taxonomy thus changes the evidence question from whether learning
``works'' to what the complete configuration establishes under a stated release
contract. The audit leaves four relationships for management research: how
effective reuse converts lifecycle investment into value; how bottlenecks and
structural change affect the allocation of binding learned and retained
authority; how error consequences, assurance, and deadline-feasible recovery
shape automated release; and how organizational rights align with information
and reversibility. The taxonomy defines these choices, and the audit identifies
the unresolved relationships. Section~\ref{sec:om_contingency} uses
operations-management theory to state their expected direction.


\section{Configuring Learning-Augmented Scheduling Systems: A Contingency Framework}
\label{sec:om_contingency}

The taxonomy and evidence audit do not support a universal ranking of
solver-led configurations, configurations with shared technical authority, and
model-led configurations. Configuration is therefore a contingency-fit problem
\citep{drazin1985fit}. The unit of choice is the complete scheduling decision
system, including its technical components and organizational roles. We define
lifecycle value as the operating benefit of released decisions less the costs
of capability acquisition, integration, decision production, assurance,
recovery, updating, and maintenance. Relative lifecycle value is the difference
between this value and that of the best feasible alternative under the same
release contract, subject to safety, legal, and business-continuity
requirements.

The taxonomy supplies the design variables in Table~\ref{tab:contingency_framework},
and the audit identifies the relationships that existing comparisons leave
unresolved. Operations-management theory supplies the directional mechanisms.
The propositions therefore state relationships to be tested rather than causal
effects already established by the reviewed studies.

\begin{table}[!htbp]
\centering
\caption{Contingency logic for configuring learning-augmented scheduling systems}
\label{tab:contingency_framework}
\begingroup
\setlength{\tabcolsep}{2.6pt}
\renewcommand{\arraystretch}{1.08}
\scriptsize
\begin{tabularx}{\textwidth}{P{2.00cm} P{3.05cm} P{3.75cm} Y}
\toprule
Management decision
& Taxonomic design variable
& Audited diagnosis
& Operations-management mechanism and proposed relationship \\
\midrule

Whether to invest
& Capability-acquisition and updating infrastructure; effective reuse volume
& Selected systems report low recurring computation and repeated use, but
acquisition, integration, updating, and maintenance lack common lifecycle
accounting.
& Fixed capability-building and integration costs are amortized over the
decisions served while the capability and supporting assets remain applicable.
\\

Allocation of technical authority
& Scope of binding learned formation; decision space retained for revision by
explicit optimization
& Reported performance is configuration-specific. Authority allocations have
not been compared under a common bottleneck and change process, and structural
transfer remains bounded.
& Binding learned formation gains value when it replaces work at a stable
formation bottleneck; retained revisable authority gains value when change
alters the formulation. \\

Scope of automated release
& Release scope; feasibility, performance, and applicability controls;
conditional recovery
& The audited records do not jointly estimate release scope, error
consequences, assurance coverage, and recovery capacity.
& Automation saves delay and release labor but adds error exposure plus
assurance and recovery costs. Its net value depends on whether acceptable
decisions can still be released before the deadline. \\

Assignment of organizational rights
& Allocation of specification, release, intervention, and lifecycle rights
& Field cases document different allocations, but do not estimate the effect
of rights--information alignment.
& Rights create value when they follow relevant information and accountability,
provided that coordination and rollback remain timely. \\

\bottomrule
\end{tabularx}

\vspace{2pt}
\parbox{\textwidth}{\scriptsize
\emph{Notes.} Taxonomic design variables describe alternative system
configurations, not stages of maturity. The audited diagnosis reports what the
review establishes; the final column states theory-informed mechanisms and
relationships proposed for testing.}
\endgroup
\end{table}

\subsection{Investment: Effective Reuse}
\label{subsec:effective_reuse}

Learning requires information and infrastructure for acquiring, integrating,
and updating capability. Selected studies report low recurring computation,
while field cases document repeated use and system integration
\citep{liang2022lenovo,kim2024photolithography,
konstantelos2022fabwide}. The relevant scale variable is not scheduling
frequency alone. \emph{Effective reuse volume} is the number of decisions served within
the joint operating coverage and useful life of the learned capability and the
assets needed to deploy it. It falls when only a subset of decisions is covered
or when operational change requires material replacement of the capability,
data pipeline, simulator, or integration.

\emph{Proposition 1 (effective reuse).} Conditional on a positive recurring
advantage over the best feasible alternative, a learning-augmented scheduling
system produces greater relative lifecycle value as effective reuse volume
increases.

\subsection{Binding Learned Authority, Revisability, and Change}
\label{subsec:authority_change_fit}

The architecture taxonomy supplies two related design variables. Binding
learned authority is the scope of schedule choices fixed by learned output that
later components cannot restore or replace. Retained revisable authority is the
decision space that an explicit schedule-forming procedure continues to
control. Both may be present in a configuration with shared technical
authority.

The proposed mechanism is information-processing fit
\citep{galbraith1974information}. Binding learned schedule formation can replace
work at a formation bottleneck when the formulation remains stable. Retained
revisable authority becomes more valuable when change enters objectives,
constraints, routing, or resource relations because those changes can be
represented and checked explicitly.

\emph{Proposition 2 (authority--change fit).} The relative lifecycle value of
expanding binding learned authority increases when schedule formation is the
binding response-time bottleneck and the scheduling formulation remains
stable. As change increasingly alters the formulation, the relative lifecycle
value of retaining decision space for revision by explicit optimization
increases.

\subsection{Automated Release: Error, Assurance, and Recovery}
\label{subsec:automated_release_fit}

The assurance taxonomy separates the claims needed to release a candidate from
recovery after the normal path fails. Field systems report action filters,
parallel scheduling, human review, and staged release
\citep{kim2024photolithography,hong2024naphtha,
konstantelos2022fabwide}. These records document the reported control
arrangements.

The mechanism is a net-value balance. Broader automated release reduces
case-by-case delay and labor. It also exposes more decisions to release error
and may require additional checking and reserved recovery capacity. A fallback
supports automation only when it can produce and recheck an acceptable
alternative before the decision deadline.

\emph{Proposition 3 (net value of automated release).} Relative to case-by-case
release, wider automated release creates greater relative lifecycle value as
error consequences fall and as assurance coverage and deadline-feasible
recovery capacity increase, after accounting for added assurance and recovery
costs.

\subsection{Governance: Rights and Reversibility}
\label{subsec:rights_alignment}

Technical authority over schedule content does not determine who may specify
the problem, approve release, intervene, or change the system. Field cases
allocate these rights differently across planners, engineers, managers, and
technical owners
\citep{liang2022lenovo,kim2024photolithography,
konstantelos2022fabwide}. The audit establishes these different arrangements,
not the performance effect of alignment. Proposition~4 is therefore derived
from authority theory \citep{aghion1997authority} and from the taxonomy's
separation of technical and organizational authority.

The proposed mechanism assigns each right to an arrangement that combines
access to the relevant information with accountability for its consequences.
Different rights may therefore have different holders. The benefit is bounded
when coordination exceeds the decision window or an intervention cannot be
reversed without substantial operational loss.

\emph{Proposition 4 (rights alignment).} Relative to an otherwise equivalent
configuration, aligning specification, release, intervention, and lifecycle
rights with relevant information and accountability increases lifecycle value
when required coordination remains timely and effective rollback is available.

\paragraph{Testing the propositions.}
\label{subsec:testing_propositions}

Tests should compare complete systems under a common evaluation contract and
estimate whether relative lifecycle value changes with the contingency named
in each proposition, rather than rank configurations by average performance.
EC.6 specifies the constructs, measurement fields, research designs, and
findings that would be inconsistent with each proposition.

\section{Conclusion}
\label{sec:conclusion}

This review treats the complete reported scheduling decision system---from
problem specification through release and recovery---as its unit of analysis.
The resulting taxonomy shifts comparison from isolated algorithms to
configurations that differ in what learning does, how schedule content is
formed, where binding technical authority resides, and how released decisions
are assured.

Applying this taxonomy to the evidence supports selected quality--computation-time trade-offs, bounded transfer of learned capability, performance under specified operating regimes, and the operational possibility of several configurations. These findings remain tied to the reported release paths and operating conditions. 
Across the audited analytical library, we did not identify a system-equivalent ranking of learning and optimization or a common basis for lifecycle economics and long-run field performance.

The contingency framework therefore changes the management question from whether learning should replace optimization to how learning, optimization, assurance, and organizational rights should be configured together. Relative lifecycle value depends on effective reuse, the fit between authority allocation and the form of change, deadline-feasible assurance and recovery, and the alignment of decision rights with information and accountability. Solver-led configurations, configurations with shared technical authority, and model-led configurations are alternative designs whose value depends on the operating setting.

\clearpage
\begingroup
\raggedbottom

\setcounter{section}{0}
\setcounter{subsection}{0}
\setcounter{table}{0}
\setcounter{figure}{0}
\renewcommand{\thesection}{EC.\arabic{section}}
\renewcommand{\thesubsection}{EC.\arabic{section}.\arabic{subsection}}
\renewcommand{\thetable}{EC.\arabic{table}}
\renewcommand{\thefigure}{EC.\arabic{figure}}
\setcounter{secnumdepth}{2}
\setcounter{tocdepth}{1}

\makeatletter
\@ifpackageloaded{hyperref}{%
  \renewcommand{\theHsection}{EC.\arabic{section}}%
  \renewcommand{\theHsubsection}{EC.\arabic{section}.\arabic{subsection}}%
  \renewcommand{\theHtable}{EC.\arabic{table}}%
  \renewcommand{\theHfigure}{EC.\arabic{figure}}%
}{}
\makeatother

\renewcommand{\arraystretch}{1.08}
\setlength{\LTpre}{3pt}
\setlength{\LTpost}{6pt}
\setlength{\tabcolsep}{3.5pt}

\makeatletter
\@ifclassloaded{poms1_V1}{%
  \RUNAUTHOR{Anonymous Authors}%
  \RUNTITLE{Electronic Companion: Learning-Augmented Scheduling}%
}{}
\makeatother

\begin{center}
{\Large\bfseries Appendix}

\vspace{7pt}
\end{center}

This electronic companion provides the review design and source-adjudicated
records supporting the main paper. EC.1 documents candidate identification,
the analytical library, and the comparison with prior reviews. EC.2 defines
the system-level codebook and resolves decisive category boundaries. EC.3
reports representative records and the 101-study common-field registry. EC.4
applies the taxonomy to the technical claims used in the evidence synthesis.
EC.5 reports the operational-maturity rules and application episodes. EC.6
links the reviewed gaps and proposed relationships to direct empirical tests.

\section{Review Approach and Analytical Evidence Base}
\label{ec:identification}

\subsection{Critical Integrative Review Design}
\label{ec:review_design}

This study is a critical integrative review. It combines structured candidate identification with a purposively assembled full-text analytical library. The
search broadens and documents discovery, including possible operational episodes. The analytical library supports taxonomy development, boundary
adjudication, representative examples, and bounded appraisal of technical, operational, and managerial claims. Search counts therefore describe the
review process; they do not estimate field composition, method frequency, or deployment prevalence.

The conceptual unit is the complete reported scheduling decision system: the
path by which an operating problem becomes a schedule that may be released,
repaired, rejected, or replaced. Bibliographic accounting uses a
\emph{substantive-work family}, which consolidates reports of the same
contribution. Materially different evaluated paths within a family are treated
as separate configurations. Evidence conclusions are stored as
claim--evidence records, and operational evidence is grouped by
system--organization--setting episode. EC.2 defines the coding units applied to
these records.

\subsection{Structured Candidate Identification}
\label{ec:search_construction}

The structured search combined broad task--resource scheduling phrases with a
general learning block. It did not use model-family, training-method,
solution-formation, assurance, authority, or deployment labels. Thirty-seven
direct-study and nine prior-review queries were executed across Semantic
Scholar, DBLP, arXiv, OpenReview, and PubMed for material available by August
26, 2026. Exact strings, field mappings, execution records, and the work-family
crosswalk are retained with the review materials.

Identifier normalization, cross-platform deduplication, and publication-version
consolidation produced 2,408 candidate work families: 2,347 from the
direct-study route and 61 from the prior-review route. The former broadens
direct-study and operational-episode discovery; the latter supports discovery
of earlier syntheses and citation tracing. A search match establishes only
candidate status. The search did not define a closed study population or the
membership of the analytical library.

\subsection{Purposive Full-Text Selection and Source-Use Criteria}
\label{ec:source_use}

The full-text library contains 148 works: 58 linked to the structured search
and 90 purposive supplements. Supplements extend conceptual variation across
interfaces, formation processes, technical-authority configurations, learning
and adaptation, assurance and recovery, claim types, and operational maturity.
They also include contrasting cases, difficult boundaries, foundational
methods, prior reviews, and operations-management theory. Search-linked and
supplemental works follow the same publication, version, cutoff, and
evidence-role rules.

Within the library, 101 analytical studies receive the common classification
reported in EC.3. Eleven additional studies receive targeted coding for
operational or domain-specific questions, and 36 sources provide contextual
methodological, benchmarking, review, application, or theoretical support.
These groups serve different analytical purposes and do not form a
representative sample.

The archival report is cited when available, while the earliest complete
public report determines availability by the cutoff. A complete stable preprint
may support an emerging result when no archival version was available. English
full text was required for source-level coding.

A work supplies \emph{direct scheduling evidence} only when it satisfies three
conditions:

\begin{itemize}[leftmargin=1.6em,itemsep=2pt,topsep=3pt]
\item \emph{P---scheduling problem}: allocation, sequencing, dispatching, or
timing under explicit temporal or resource relations and a stated scheduling
objective;
\item \emph{L---reusable learning}: a learned artifact persists beyond the run
that created it and is evaluated on a later episode, independent instance, or
subsequent decision epoch; and
\item \emph{D---decision effect}: the learned output affects, and is evaluated
through, a scheduling decision or an inference-time scheduling computation.
\end{itemize}

\emph{Direct} evidence satisfies P--L--D. \emph{Enabling} evidence supplies a
mechanism or benchmark used to interpret scheduling systems. \emph{Adjacent}
evidence supports only a specified comparison or research direction from a
related optimization domain. Coding may be full, targeted, or contextual
according to assigned use; operational maturity is assessed separately.

\begin{table}[!htbp]
\centering
\caption{Review records, accounting, and analytical use}
\label{tab:ec_evidence_sets}
\small
\begin{tabularx}{\textwidth}{L{3.25cm}L{3.10cm}Y}
\toprule
Record & Accounting & Analytical use \\
\midrule
Structured candidate record
& 2,408 work families: 2,347 direct-route and 61 prior-review-route
& Discovery and source tracing; not a full-text inclusion set or a basis for
frequency estimates. \\
\addlinespace[2pt]
Full-text analytical library
& 148 works: 58 search-linked and 90 purposive supplements
& Taxonomy development, boundary analysis, representative cases, and bounded
claim appraisal. \\
\addlinespace[2pt]
Common-field registry
& 101 analytical studies
& Comparison through common taxonomy fields; materially different paths and
unresolved classifications remain visible. \\
\addlinespace[2pt]
Targeted and contextual records
& 11 targeted studies and 36 contextual sources
& Operational or domain-specific adjudication and support for specified
methodological, benchmarking, review, and theoretical claims. \\
\bottomrule
\end{tabularx}
\end{table}

\subsection{Targeted Equipment-Specific Records}
\label{ec:targeted_equipment}

The 11 targeted studies comprise the five Level A records reported in EC.5 and
the six equipment-specific scheduling studies summarized below. The latter
extend coverage of cluster-tool, robotic-cell, and robotic-flow-shop
scheduling but remain outside the 101-study common-field registry.

Each study provides direct scheduling evidence for its reported equipment
model. Computational evaluation, domain detail, or parameter calibration does
not by itself establish Level A operating use or Level B1 bounded physical
execution.

\begin{center}
\begingroup
\setlength{\tabcolsep}{3pt}
\renewcommand{\arraystretch}{1.08}
\scriptsize
\begin{tabularx}{\textwidth}{P{3.25cm} P{4.60cm} Y}
\toprule
Study and setting
& Reported learning configuration
& Evidence boundary \\
\midrule

\citet{kim2025noncyclic}; cluster tools
& Deep reinforcement learning (DRL) for noncyclic scheduling across multiple tool
and wafer-flow configurations
& Computational evaluation on generated instances; the report does not
document live-tool release or sustained operating use. \\

\addlinespace[2pt]
\citet{kim2025concurrent}; cluster tools
& Deep Q-learning with adaptive look-ahead search for concurrent processing
& Direct computational scheduling evidence; production release and field
duration are not documented. \\

\addlinespace[2pt]
\citet{kim2024roboticcell}; robotic cells
& Deep reinforcement learning with look-ahead search for single-gripper
robotic-cell scheduling
& Computational evaluation, including processing-time variation; physical or
sustained field operation is not documented. \\

\addlinespace[2pt]
\citet{kim2022dualgripper}; dual-gripper robotic cells
& Reinforcement learning (RL) for scheduling dual-gripper cells
& Experimental parameters are partly calibrated from semiconductor data, but
schedule evaluation remains computational and does not document live-cell
integration. \\

\addlinespace[2pt]
\citet{lee2022roboticflowshop}; robotic flow shops
& Reinforcement learning under modeled processing-time variation
& Direct computational scheduling evidence; field execution, organizational
roles, and operating duration are not documented. \\

\addlinespace[2pt]
\citet{kim2023lookahead}; robotic flow shops
& Deep reinforcement learning with look-ahead search for noncyclic
dual-gripper flow shops
& Computational evaluation includes bounded train--test changes but does not
establish production deployment or plant-level performance. \\

\bottomrule
\end{tabularx}
\endgroup
\end{center}

\subsection{Positioning Relative to Prior Reviews}
\label{ec:prior_review_positioning}

Table~\ref{tab:ec_prior_reviews} compares the present organizing questions
with ten prior reviews in the analytical library. The table identifies
differences in focus rather than ranking review quality. A dimension is
\emph{core} when it organizes a research question, dedicated section, or
repeated study-level comparison; \emph{partial} denotes explicit but bounded
treatment; and \emph{outside} means that the dimension is not a focal concern.

{\scriptsize
\setlength{\tabcolsep}{2.4pt}
\begin{longtable}{L{2.50cm}L{5.20cm}ccccc}
\caption{Positioning relative to prior reviews}
\label{tab:ec_prior_reviews}\\
\toprule
Review & Principal scope and organizing unit & IA & LA & AR & EO & GD \\
\midrule
\endfirsthead
\multicolumn{7}{l}{\small\itshape Table \thetable\ continued}\\
\toprule
Review & Principal scope and organizing unit & IA & LA & AR & EO & GD \\
\midrule
\endhead
\bottomrule
\endfoot
\citet{bengio2021methodological}
& Generic combinatorial optimization (CO); learning source and machine learning (ML)--optimization structure
& C & C & P & P & O \\
\citet{kotary2021survey}
& Constrained optimization; ML-augmented and end-to-end integration patterns
& C & P & P & P & O \\
\citet{cappart2023combinatorial}
& Graph-based combinatorial optimization; primal, dual, and algorithmic-reasoning roles
& C & C & P & P & O \\
\citet{smit2025gnnsurvey}
& Job-shop and flow-shop scheduling; graph representation, task, and training method
& P & C & P & P & O \\
\citet{lv2025drlreview}
& Job-shop and flexible-job-shop scheduling; state, action, network, and RL design
& P & C & O & P & O \\
\citet{khadivi2025drlreview}
& Machine scheduling; DRL component and machine environment
& P & C & P & P & O \\
\citet{zhang2025rlreview}
& Job-shop scheduling; RL algorithm, agent structure, and policy update
& P & C & O & P & O \\
\citet{goppert2024online}
& Online manufacturing scheduling; problem, ML type, and objective
& O & C & O & P & O \\
\citet{zhang2024resilientreview}
& Dynamic manufacturing; application setting and resilience or sustainability indicator
& P & C & P & P & P \\
\citet{stockermann2025deploymentreview}
& Semiconductor dispatching; validation, deployment challenge, and implementation pipeline
& P & C & C & C & P \\
\end{longtable}
}

\noindent\footnotesize\textit{Notes.}
C = core; P = partial; O = not a focal analytical dimension. IA = inference
architecture; LA = learning and adaptation; AR = assurance and release; EO =
evidence appraisal and operational maturity; GD = operations-management
governance and decision rights.\normalsize

The reviewed syntheses provide complementary foundations in optimization,
scheduling, and deployment. The present review connects these concerns at the
level of a released scheduling decision system, including downstream work,
normal-path technical authority, assurance and recovery, evidence strength,
and organizational decision rights.

\section{System-Level Codebook and Boundary Adjudication}
\label{ec:codebook}

\subsection{Coding Structure and Source Standard}
\label{ec:relational_structure}

The codebook distinguishes a substantive work from its reports and then
separates materially different evaluated release paths:

\begin{quote}
\small
substantive-work family $\rightarrow$ adjudicated report $\rightarrow$
evaluated configuration $\rightarrow$ learned link, formation stage, update
phase, and assurance record $\rightarrow$ technical claim or application
episode.
\end{quote}

One work may therefore supply several analytical records without being counted
several times bibliographically. Every coded fact is tied to the report version
and source passages used for adjudication. \texttt{NR} denotes an applicable
fact not reported in the checked source; \texttt{NA} denotes a structurally
inapplicable field; and \texttt{UNRESOLVED} denotes a field for which the source
does not support a defensible classification. Model names and author
terminology do not determine a code.

\begingroup
\footnotesize
\begin{longtable}{L{2.75cm}L{5.25cm}L{6.60cm}}
\caption{System-level codebook}
\label{tab:ec_codebook}\\
\toprule
Record and coding unit & Question and recorded values & Boundary rule \\
\midrule
\endfirsthead
\multicolumn{3}{l}{\small\itshape Table \thetable\ continued}\\
\toprule
Record and coding unit & Question and recorded values & Boundary rule \\
\midrule
\endhead
\bottomrule
\endfoot

Operating boundary; configuration and claim
& Under what scheduling structure and information regime was the result
tested? Resource, precedence, routing, setup, and capacity structure;
objective and preference regime; information timing; scale; tested change.
& These conditions bound every code and claim. Same-generator size testing is
not structural transfer, and a sequential Markov decision process is not
dynamic scheduling unless information arrives during operation. \\

Learning interface; learned output--consumer link
& What does the first downstream consumer use the output to do? SUB surrogate
substitution; PROP candidate proposal; INIT initialization; SEL component
selection or configuration; SC-O order-preserving search control; SC-C
coverage-changing search control; FORM direct schedule formation.
& Code downstream use, not output shape or model name. Split materially
different consumers. A priority vector is FORM under fixed decoding but SC-O
or SC-C when retained search remains responsible for construction. \\

Solution formation; learned stage within a normal release path
& Does learning form inherited schedule content, and how? SG0 no learned
release-level formation; SG1 simultaneous one-shot; SG2 successive conditional
construction; SG3 learned revision of an executable candidate; SG4 repeated
global update of a nonexecutable latent or relaxed state.
& SG1--SG4 classify a stage; one path may compose stages. An overwritable warm
start is SG0. Repeated rolling-horizon releases do not create SG2, and
training-time self-improvement is not SG3. \\

Normal-path technical authority; complete normal release path
& Which components exert binding control over releasable schedule content?
Solver-led (SL); shared technical authority (SH); model-led (MDL).
& SL has no binding learned contribution. SH combines binding learned content
with independent retained formation. MDL has binding learned content followed
only by fixed decoding, propagation, checking, or veto. Recovery and
organizational rights are separate. \\

Learning and adaptation; update term or phase
& What task information directly defines the update? Independent-reference
target (IR); model-generated selected target (MG); differentiable problem loss
(DP); return or fitness (RF). Acquisition, timing, updated object, and
concurrent or ordered composition are modifiers.
& Oracle targets are IR; objective-ranked pseudo-labels are MG; gradients from
problem losses are DP; rewards, advantages, Bellman targets, and black-box
fitness are RF. Pretext relations defined independently of the current learner
are auxiliary reference signals; generic regularizers are recorded separately. \\

Decision assurance; claim--candidate/path record
& Which outcome is required, by what control, under what evidence and scope,
and with what failure response? F feasibility, P performance, and A
applicability are three claims supporting release; R is the conditional
recovery path after a failed check or missed deadline.
& F is not a performance certificate, and empirical confidence is not a hard
feasibility check. R requires a trigger, available alternative-producing
capacity, rechecking, and deadline evidence. A recovery candidate receives
separate F/P/A records. \\

Organizational decision rights; configuration or episode
& Who holds specification, release, intervention, and lifecycle rights?
Reported holder, formal right or effective control, relevant information,
accountability, coordination, and rollback.
& Rights are not inference-architecture categories and cannot be inferred from
technical authority or human presence. Unreported holders remain \texttt{NR}. \\
\end{longtable}
\endgroup

Technical claim fields are applied in EC.4, operational maturity in EC.5, and
the proposition constructs in EC.6; they are not additional taxonomy
categories.

\subsection{Coding Order and Decision Rules}
\label{ec:coding_order}

Coding follows the reported pipeline. After the work family, report,
scheduling problem, and evidence role are established, materially different
normal release paths are separated. Each learned output is traced to its first
consumer; SG1--SG4 are assigned only to learned content inherited by that path;
and technical authority is assigned after all binding contributions have been
identified. Update, assurance, claim, episode, and organizational-rights
records are then linked to the relevant configuration.

A \emph{binding schedule effect} fixes or removes a scheduling choice that no
later component on the same path has both the structural right and remaining
decision space to restore or replace. Thus an overwritable warm start differs
from learned fixing, and a checker that can only accept or veto a candidate
supplies assurance rather than shared authority. Conditional fallback does not
change the normal-path authority code.

Formation follows the state being changed rather than the number of model
calls. SG1 creates learned content simultaneously; SG2 conditions later learned
choices on earlier commitments; SG3 revises a complete executable candidate;
and SG4 repeatedly updates a global nonexecutable representation. Learning is
coded one separable update term or phase at a time, with concurrent terms and
ordered phases retained explicitly.

The registry was developed through a sequential source audit rather than
duplicate independent coding. The first pass covered all 101 common-field
studies. A targeted second pass rechecked 59 works used for representative
configurations, difficult boundaries, evidence conclusions, or managerial
claims. The technical-authority fields for \citet{liu2024integrating} and
\citet{tassel2023endtoend} remain \texttt{UNRESOLVED}; the review reports
neither an intercoder-agreement statistic nor inferred values for unreported
fields.

\subsection{Boundary Adjudications}
\label{ec:boundary_cases}

Table~\ref{tab:ec_boundary_cases} retains six cases in which the reported
downstream path, formation state, or update route changes the classification.

\begingroup
\footnotesize
\begin{longtable}{L{3.25cm}L{5.85cm}L{5.70cm}}
\caption{Decisive boundary adjudications}
\label{tab:ec_boundary_cases}\\
\toprule
Configuration & Decisive pipeline fact & Adjudicated record and boundary \\
\midrule
\endfirsthead
\multicolumn{3}{l}{\small\itshape Table \thetable\ continued}\\
\toprule
Configuration & Decisive pipeline fact & Adjudicated record and boundary \\
\midrule
\endhead
\bottomrule
\endfoot

One-shot priorities in \citet{jeon2023neural}; prediction with  completion by linear programming (LP)
in \citet{kotary2022fast}
& Jeon applies fixed list decoding. Kotary preserves learned machine order but
uses normal-path LP completion to determine timing.
& Both contain SG1 formation. Jeon is FORM/SG1/MDL; Kotary is FORM/SG1/SH.
Fixed decoding does not create shared authority, whereas independent retained
completion does. \\

Three paths in \citet{omalley2023hybrid}
& Stochastic mixed-integer program (MIP) may overwrite the learned warm start; a model-free policy
forms commitments; a lookahead variant restricts actions before scenario
evaluation selects a commitment.
& Hybrid: INIT/SG0/SL. Model-free policy: FORM/SG2/MDL. Lookahead:
SC-C/SG0/SH. Each evaluated path is a separate configuration. \\

Predicted precedences in \citet{verhaeghe2024precedences}
& The same predictions are either added as constraints or used only to order
constraint programming (CP) search.
& Constraint addition: FORM/SG1/SH because choices are removed. Ordering:
SC-O/SG0/SL because the original decision space remains available. \\

Constructive policies in \citet{zhang2020dispatch} and
\citet{corsini2024selflabeling}
& Both form schedules through successive learned choices, but one updates from
sampled returns and the other fits selected model-generated candidates.
& Both are FORM/SG2/MDL at inference; their update records are RF and MG,
respectively. Formation does not determine the learning mechanism. \\

Diffusion-based formation in \citet{li2026goaldiffusion}, with
\citet{sun2023difusco} as adjacent mechanism evidence
& Repeated denoising changes a full nonexecutable representation before fixed
decoding rather than extending a partial schedule or revising an incumbent.
& FORM/SG4/MDL for the direct scheduling configuration. The scheduling result
was an emerging preprint at the cutoff; routing results remain adjacent. \\

Per-instance differentiation in \citet{liu2024differentiable}; reusable
learning with adaptation in \citet{liu2025integrated}
& The first optimizes new decision variables for each instance and retains no
learned artifact. The second trains a reusable model and then adapts a
parameter copy for the test instance.
& First: DP mechanism evidence with reusable-learning configuration
\texttt{NA}. Second: ordered IR $\rightarrow$ DP; adaptation belongs in the
decision-work and cost record. \\
\end{longtable}
\endgroup

\section{Taxonomy Records and the 101-Study Common-Field Registry}
\label{ec:taxonomy_registry}

EC.3 applies the codebook to source-adjudicated literature records. The first
table is a configuration crosswalk covering all seven interface codes,
SG0--SG4, and the three technical-authority configurations. A second table
records representative update phases. The compact registry then retains the
common taxonomy fields for all 101 analytical studies. Claims and operational
episodes are assessed in EC.4 and EC.5 rather than inferred from code counts.

\subsection{Representative Configuration and Interface Records}
\label{ec:interface_records}
\label{ec:configuration_records}

Table~\ref{tab:ec_configuration_records} reports materially different normal
release paths rather than assigning one headline category to a study. F, P,
A, and R tags identify the assurance issue represented in the concise record;
they do not imply that every assurance claim was established.

\begingroup
\scriptsize
\setlength{\tabcolsep}{2pt}
\renewcommand{\arraystretch}{1.07}
\begin{longtable}{L{3.25cm}L{1.45cm}L{2.30cm}L{1.55cm}L{6.15cm}}
\caption{Source-adjudicated representative configurations}
\label{tab:ec_configuration_records}\\
\toprule
Configuration & Interface & Formation; authority & Update & Assurance and decisive boundary \\
\midrule
\endfirsthead
\multicolumn{5}{l}{\small\itshape Table \thetable\ continued}\\
\toprule
Configuration & Interface & Formation; authority & Update & Assurance and decisive boundary \\
\midrule
\endhead
\bottomrule
\endfoot

Surrogate-assisted decomposition in \citet{bouska2023deep}
& SUB & SG2; SH & IR
& F/P: the retained recursive heuristic makes the scheduling commitments;
A/R: \texttt{NR}. \\

Learned column proposal in \citet{hijazi2026improving}
& PROP & SG0; SL & IR
& F/P: the restricted master and exact pricing evaluate optional columns and
establish only the reported master-problem claim; A/R: \texttt{NR}. \\

Learned warm-start hybrid in \citet{omalley2023hybrid}
& INIT & SG0; SL & RF
& F/P: retained stochastic MIP may overwrite the learned schedule and supplies
the reported status or gap; A: finite-scenario formulation; R: \texttt{NR}. \\

CP-solver selection in \citet{muller2022algorithm}
& SEL & SG0; SL & IR
& F/P: the selected CP solver forms the schedule and supplies bound or proof
status when achieved; A/R: \texttt{NR}. \\

Learned ordering within CP search in \citet{sun2024enhancing}
& SC-O & SG0; SL & IR
& F/P: the original CP model and branches remain available; solver status, not
the prediction, determines the guarantee; A/R: \texttt{NR}. \\

Policy-guided lookahead in \citet{omalley2023hybrid}
& SC-C & SG0; SH & RF
& P/A: scenario evaluation selects among retained actions after learning
restricts coverage; F/R: no general release-level claim. \\

Neural heuristic measures within ant colony optimization (ACO) in \citet{ye2023deepaco}
& SC-O & SG0; SL & RF
& F/P: retained ants, pheromone updates, and selection form and evaluate
schedules; learned measures do not create SG2; A/R: \texttt{NR}. \\

Predicted precedences as constraints in \citet{verhaeghe2024precedences}
& FORM & SG1; SH & IR
& F/P: CP acts on the augmented model and cannot restore removed alternatives;
any certificate is limited to that model; A/R: \texttt{NR}. \\

One-shot priority formation in \citet{jeon2023neural}
& FORM & SG1; MDL & RF
& F: fixed list decoding enforces represented rules; P: empirical objective
evaluation; A/R: \texttt{NR}. \\

Prediction with LP completion in \citet{kotary2022fast}
& FORM & SG1; SH & IR+DP
& F: normal-path LP completion for a consistent learned order; R: a greedy
branch was defined but not triggered; P/A: no general certificate. \\

Learning-guided rolling-horizon optimization in \citet{li2025rolling}
& FORM & SG1; SH & IR
& F: CP-SAT completes the restricted release problem after learned fixing;
P/A: bounded rolling-horizon evaluation; R: \texttt{NR}. \\

Constructive dispatching in \citet{zhang2020dispatch}
& FORM & SG2; MDL & RF
& F: eligible-action construction covers modeled constraints; P: empirical
comparison; A/R: \texttt{NR}. \\

Learned schedule improvement in \citet{zhang2024improvement}
& FORM & SG3; SH & RF
& F/P: a retained initializer supplies an executable schedule and learned N5
moves revise it; total work includes initialization and all iterations. \\

Diffusion-based schedule formation in \citet{li2026goaldiffusion}
& FORM & SG4; MDL & IR
& F: fixed topological and earliest-feasible decoding for represented rules;
P/A: emerging direct evidence at the cutoff; R: \texttt{NR}. \\

Learned-candidate release in \citet{kim2024photolithography}
& FORM & SG2; MDL & RF
& P: fixed key performance indicator (KPI) comparison selects the learned candidate on this branch; A:
documented operating use; F/R: source-bounded or \texttt{NR}. \\

Rule-based-candidate release in \citet{kim2024photolithography}
& PROP & SG0; SL & RF for rejected learned candidate
& P/A: the learned candidate is compared but rejected; on this branch it has
no binding effect. This is a distinct normal release path, not recovery. \\
\end{longtable}
\endgroup

\subsection{Representative Update Records}
\label{ec:update_records}

Table~\ref{tab:ec_update_records} assigns a mechanism to one separable update
term or phase; acquisition, timing, updated object, and composition remain
modifiers.

\begin{table}[H]
\centering
\caption{Source-adjudicated learning and adaptation records}
\label{tab:ec_update_records}
\begingroup
\setlength{\tabcolsep}{3pt}
\renewcommand{\arraystretch}{1.08}
\scriptsize
\begin{tabularx}{\textwidth}{P{4.35cm} P{3.35cm} Y}
\toprule
Study and update phase & Mechanism & Decisive update route \\
\midrule
Supervised rule learning in \citet{ingimundardottir2011supervised}
& IR & Optimal schedules supply fixed targets for a reusable dispatching rule. \\
Self-labeling in \citet{corsini2024selflabeling}
& MG & Model samples are evaluated and the selected sample becomes a fixed
pseudo-label. \\
Concurrent target and problem-loss learning in \citet{kotary2022fast}
& Concurrent IR+DP & Fixed schedule targets and differentiable constraint
violations enter one training objective as separable terms. \\
Supervised initialization followed by adaptation in \citet{liu2025integrated}
& Ordered IR $\rightarrow$ DP & A reusable model is trained first; a parameter
copy is then adapted for the current instance through a differentiable
surrogate. \\
Offline return-based learning in \citet{vanremmerden2025offline}
& RF from fixed experience & Rewards and Bellman targets come from a fixed set
of CP-generated trajectories. \\
Simulator-interactive learning in \citet{zhang2020dispatch}
& RF from simulated interaction & proximal policy optimization (PPO) uses sampled scheduling episodes and a
shaped makespan-related reward. \\
\bottomrule
\end{tabularx}
\endgroup
\end{table}

Objective-based selection followed by fitting a fixed candidate is MG rather
than RF. A differentiable loss applied only to instance-specific decision
variables is retained as a DP mechanism boundary rather than counted as a
reusable learning-augmented scheduling configuration.

\subsection{The 101-Study Common-Field Classification Registry}
\label{ec:complete_study_registry}

Table~\ref{tab:ec_study_atlas} retains all 101 analytical studies: 88 direct,
12 enabling, and one adjacent under the P--L--D rule. Each row reports the
study, evidence role and setting, inference architecture, update mechanism, and
a concise assurance or claim boundary. The tags F, P, A, and R identify the
release claim or conditional recovery issue explicitly represented in the
row; they are not prevalence counts or certificates. Detailed technical claims
and operating episodes remain in EC.4 and EC.5.

\begingroup
\scriptsize
\setlength{\tabcolsep}{2pt}
\renewcommand{\arraystretch}{1.04}
\begin{longtable}{L{3.85cm}L{3.30cm}L{1.65cm}L{5.65cm}}
\caption{The 101-study common-field classification registry}
\label{tab:ec_study_atlas}
\label{tab:ec_extended_registry}\\
\toprule
Study, evidence role, and setting & Interface; formation; authority & Update & Concise assurance or necessary boundary \\
\midrule
\endfirsthead
\multicolumn{4}{l}{\small\itshape Table \thetable\ continued}\\
\toprule
Study, evidence role, and setting & Interface; formation; authority & Update & Concise assurance or necessary boundary \\
\midrule
\endhead
\bottomrule
\endfoot

\citet{vinyals2015pointer} \newline \emph{Enabling; Generic variable-size output tasks}
& FORM; SG2; MDL
& IR
& \textbf{F:} Pointer decoding establishes variable-size autoregressive output; it supplies no scheduling-specific feasibility evidence. \\
\addlinespace[1pt]

\citet{bello2016neural} \newline \emph{Enabling; two-dimensional Euclidean traveling-salesperson problem (TSP) and 0--1 knapsack}
& FORM; SG2; MDL
& RF
& \textbf{P:} All samples, objective calls, and policy updates belong to test-time work. \\
\addlinespace[1pt]

\citet{kool2019attention} \newline \emph{Enabling; Routing problems}
& FORM; SG2; MDL
& RF
& \textbf{F:} Problem-specific masking and decoding; no direct scheduling experiment. \\
\addlinespace[1pt]

\citet{kwon2020pomo} \newline \emph{Enabling; Routing and other CO problems}
& FORM; SG2; MDL
& RF
& \textbf{Boundary:} Each start is autoregressive; deterministic best-of-model-samples evaluation does not add independent authority. \\
\addlinespace[1pt]

\citet{kwon2021matnet} \newline \emph{Direct; Flexible flow-shop scheduling (FFSP)}
& FORM; SG2; MDL
& RF
& \textbf{Boundary:} Relation-matrix encoding and problem-specific constructive decoding provide the direct NCO-to-scheduling bridge. \\
\addlinespace[1pt]

\citet{zhang2020dispatch} \newline \emph{Direct; Static deterministic job-shop scheduling (JSSP)}
& FORM; SG2; MDL
& RF
& \textbf{F/P:} Admissible-operation construction covers the modeled precedence and unary resource conflicts; it is not an optimality certificate. \\
\addlinespace[1pt]

\citet{park2021learning} \newline \emph{Direct; JSSP}
& FORM; SG2; MDL
& RF
& \textbf{F:} State-conditioned graph dispatch and action masking target the modeled basic JSSP constraints. \\
\addlinespace[1pt]

\citet{kotary2022fast} \newline \emph{Direct; JSSP}
& FORM; SG1; SH
& IR+DP
& \textbf{F/P/R:} LP completion covers consistent learned orders; greedy recovery was defined but not triggered. \\
\addlinespace[1pt]

\citet{jeon2023neural} \newline \emph{Direct; directed acyclic graph (DAG) scheduling}
& FORM; SG1; MDL
& RF
& \textbf{F/P:} One simultaneous priority sample is converted by a fixed list scheduler; the decoder enforces the represented scheduling rules but supplies no optimality certificate. \\
\addlinespace[1pt]

\citet{song2023flexible} \newline \emph{Direct; Flexible job-shop scheduling (FJSP)}
& FORM; SG2; MDL
& RF
& \textbf{F/P:} Sequential eligible operation--machine selection; assurance is bounded by the constraints represented in the action construction. \\
\addlinespace[1pt]

\citet{wang2024dualattention} \newline \emph{Direct; FJSP}
& FORM; SG2; MDL
& RF
& \textbf{F/P/A:} Dynamic operation--machine encoding and structured actions support modeled feasibility; they do not certify global optimality. \\
\addlinespace[1pt]

\citet{zhang2024improvement} \newline \emph{Direct; JSSP}
& FORM; SG3; SH
& RF
& \textbf{P:} A nonlearned dispatching rule forms the initial complete schedule and the learned policy applies repeated N5 moves; total cost includes initialization and all improvement iterations. \\
\addlinespace[1pt]

\citet{corsini2024selflabeling} \newline \emph{Direct; JSSP}
& FORM; SG2; MDL
& MG
& \textbf{A:} Accepted rollouts refresh training targets, while deployed inference remains constructive; training-time self-labeling is not a SG3 stage. \\
\addlinespace[1pt]

\citet{pirnay2024selfimprovement} \newline \emph{Direct; TSP, capacitated vehicle routing problem (CVRP), and JSSP}
& FORM; SG2; MDL
& MG
& \textbf{P:} Gumbeldore uses objective feedback between complete-sample rounds and returns the best sample; it does not revise an incumbent. \\
\addlinespace[1pt]

\citet{choi2025multipolicy} \newline \emph{Direct; Multiple job-scheduling variants}
& FORM; SG2; MDL
& MG
& \textbf{P/A:} Multiple latent-conditioned policies construct schedules; self-evolution is training-time and evaluated transfer remains bounded to the reported families. \\
\addlinespace[1pt]

\citet{wang2025mistar} \newline \emph{Direct; FJSP}
& INIT, PROP; SG2--SG3; MDL / nonlearned-initializer variant: PROP; SG3; SH
& RF
& \textbf{P:} Authority differs with the learned or nonlearned initializer; runtime includes initialization and all improvement iterations. \\
\addlinespace[1pt]

\citet{xiao2026resched} \newline \emph{Direct; FJSP, JSSP, and FFSP}
& FORM; SG2; MDL
& RF
& \textbf{F/A:} Models are trained by problem, so results establish framework reuse rather than unrestricted one-policy transfer. \\
\addlinespace[1pt]

\citet{luttmann2026multiaction} \newline \emph{Direct; FJSP, FFSP, and an adjacent routing task}
& FORM; SG2; MDL
& MG
& \textbf{A:} Joint assignments are made at each of several construction steps; training self-improvement does not create deployed SG3. \\
\addlinespace[1pt]

\citet{smit2025stochastic} \newline \emph{Direct; Stochastic FJSP}
& FORM; SG2; MDL
& RF
& \textbf{F/A:} Scenario processing conditions a constructive base policy; validity follows the decoder, while robustness remains conditional on sampled scenario and distribution assumptions. \\
\addlinespace[1pt]

\citet{smit2026multiobjective} \newline \emph{Direct; Multi-objective FJSP}
& FORM; SG2; MDL
& RF
& \textbf{F/P:} Preference-conditioned constructive decoding supplies feasible candidates; Pareto quality and coverage are empirical rather than certified. \\
\addlinespace[1pt]

\citet{liu2025integrated} \newline \emph{Direct; Single-machine scheduling with non-decreasing min-sum objectives, with or without release times} \newline \emph{Publication status: preprint at the cutoff}
& FORM; SG1; MDL
& IR $\rightarrow$ DP
& \textbf{P/A:} A continuous descriptor is mapped by a fixed discrete construction interface; adaptation cost is part of inference and the post-adaptation generation stage remains SG1. \\
\addlinespace[1pt]

\citet{bouska2023deep} \newline \emph{Direct; Single-machine total tardiness}
& SUB; SG2; SH
& IR
& \textbf{Boundary:} Predicted criterion values replace calculations within successive Lawler/symmetric-decomposition commitments; the nonlearned heuristic assembles the released schedule. \\
\addlinespace[1pt]

\citet{liu2024integrating} \newline \emph{Direct; JSSP}
& PROP; SG0; authority \textnormal{\texttt{UNRESOLVED}}
& IR
& \textbf{F/R:} Accepted proposals use surrogate Lagrangian relaxation coordination; normal-path authority remains \texttt{UNRESOLVED}, and branch-and-cut is conditional recovery. \\
\addlinespace[1pt]

\citet{hijazi2026improving} \newline \emph{Direct; Unrelated parallel-machine scheduling with total weighted completion time}
& PROP; SG0; SL
& IR
& \textbf{P/A:} Transformer-generated columns are internal objects; the restricted master evaluates them and exact dynamic-programming pricing at termination restores the stated column-generation LP claim, not integer global optimality by itself. \\
\addlinespace[1pt]

\citet{sun2024enhancing} \newline \emph{Direct; JSSP CP benchmarks}
& SC-O; SG0; SL
& IR
& \textbf{F/P:} Predicted solution information orders CP decisions while complete search retains feasibility and the declared search space; exactness requires certified termination. \\
\addlinespace[1pt]

\citet{heinz2025reinforcement} \newline \emph{Direct only for the ratings-carryover configuration; Job-shop and resource-constrained project-scheduling search}
& SC-O; SG0; SL
& RF
& \textbf{P:} Certification remains conditional on complete solver termination, and only the cross-run carryover configuration supplies direct reusable-learning evidence. \\
\addlinespace[1pt]

\citet{strassl2022instance} \newline \emph{Direct; JSSP}
& SEL; SG0; SL
& IR
& \textbf{P:} The selected nonlearned algorithm forms the schedule; error consequences are misrouting and oracle regret rather than learned schedule generation. \\
\addlinespace[1pt]

\citet{muller2022algorithm} \newline \emph{Direct; FJSP}
& SEL; SG0; SL
& IR
& \textbf{F:} The selected CP solver retains schedule formation, feasibility enforcement, and search responsibility. \\
\addlinespace[1pt]

\citet{omalley2023hybrid} \newline \emph{Direct; Stochastic unit-commitment scheduling}
& Hybrid: INIT; SG0; SL / model-free RL: FORM; SG2; MDL / lookahead RL: SC-C; SG0; SH
& RF
& \textbf{A:} Warm-start, model-free, and lookahead variants are distinct paths; every sample, simulation, and retained solve counts. \\
\addlinespace[1pt]

\citet{tassel2023endtoend} \newline \emph{Direct; JSSP}
& FORM; SG2; authority \textnormal{\texttt{UNRESOLVED}}
& IR+RF
& \textbf{F:} CP propagation supports modeled validity; authority remains \texttt{UNRESOLVED} because an additional selector's role is unclear. \\
\addlinespace[1pt]

\citet{ingimundardottir2011supervised} \newline \emph{Direct; JSSP dispatching}
& FORM; SG2; MDL
& IR
& \textbf{P:} The learned linear priority rule is applied successively to partial shop states; dispatching semantics form the schedule and performance inherits label and feature coverage. \\
\addlinespace[1pt]

\citet{ingimundardottir2018discovering} \newline \emph{Direct; Job-shop dispatching}
& FORM; SG2; MDL
& IR
& \textbf{A:} A learned dispatching rule acts repeatedly during construction; iterative data collection mitigates but does not eliminate state-distribution shift. \\
\addlinespace[1pt]

\citet{parmentier2023structured} \newline \emph{Direct; Single machine with release times and total completion time}
& SUB; SG0; SH
& IR
& \textbf{P:} A predictor transforms the hard instance into an easier instance that a nonlearned algorithm solves optimally before transposition to the original problem; the learned output is not itself a release-level schedule. \\
\addlinespace[1pt]

\citet{li2024imitation} \newline \emph{Direct; Permutation flow-shop scheduling (PFSP)}
& FORM; SG2; MDL
& IR
& \textbf{F/P:} Masking prevents reselection, so decoding yields a valid permutation; the improved Nawaz--Enscore--Ham (NEH) heuristic supplies offline labels, while inference has no quality or optimality certificate. \\
\addlinespace[1pt]

\citet{liu2024differentiable} \newline \emph{Enabling mechanism boundary; Latency-constrained minimum-resource scheduling of dataflow DAGs}
& Reusable-learning architecture: \textnormal{\texttt{NA}}
& DP applied to instance-specific decision variables
& \textbf{Boundary:} No reusable artifact is retained; the work supplies DP mechanism evidence with reusable architecture \texttt{NA}. \\
\addlinespace[1pt]

\citet{vanremmerden2025offline} \newline \emph{Direct; JSSP dispatching}
& FORM; SG2; MDL
& RF
& \textbf{F:} Maskable actions support modeled validity; CP is used for data generation, not inference-time completion, and behavior outside logged support remains a risk. \\
\addlinespace[1pt]

\citet{luo2020dynamic} \newline \emph{Direct; Dynamic FJSP with new-job arrivals}
& SEL; SG2; MDL
& RF
& \textbf{F/A:} Repeated rule selection constructs the executed schedule; the simulator enforces modeled transitions and dynamic validity depends on event fidelity. \\
\addlinespace[1pt]

\citet{luo2021multiobjective} \newline \emph{Direct; Dynamic multi-objective FJSP}
& SEL; SG2; MDL
& RF
& \textbf{F/P:} State-conditioned rule selection is repeated at scheduling points; feasibility and trade-offs are bounded by the rule library, simulator, and specified reward. \\
\addlinespace[1pt]

\citet{yuan2025multiobjective} \newline \emph{Direct; Dynamic multi-objective FJSP}
& FORM; SG2; MDL
& RF
& \textbf{F/P/A:} Performance remains conditional on the simulated arrivals and scalarized reward. \\
\addlinespace[1pt]

\citet{zhang2026uncertainty} \newline \emph{Direct; JSSP under structural uncertainty}
& FORM; SG2; MDL
& RF
& \textbf{F/A:} The actor dispatches operations under evolving route realization; uncertainty perception and admissible actions are evaluated only for the reported structural uncertainty process. \\
\addlinespace[1pt]

\citet{yang2026dynamic} \newline \emph{Direct; Dynamic FJSP with transportation and priority constraints}
& SEL; SG2; MDL
& RF
& \textbf{F/P/A:} Assurance is bounded by the rule library and modeled arrival/urgent-order process. \\
\addlinespace[1pt]

\citet{drakulic2025goal} \newline \emph{Direct; Multi-task CO including JSSP and unrelated-machine scheduling, with transfer to open shop}
& FORM; SG2; MDL
& IR+MG
& \textbf{P/A:} Reported fine-tuning supports bounded adaptation-based transfer, not unrestricted zero-shot scheduling transfer. \\
\addlinespace[1pt]

\citet{hottung2025polynet} \newline \emph{Direct; Multiple CO tasks including FFSP}
& FORM; SG2; MDL
& RF
& \textbf{P:} Strategy-conditioned policies construct solutions with problem-specific decoders; diversity within a problem is not cross-problem certification. \\
\addlinespace[1pt]

\citet{ye2023deepaco} \newline \emph{Direct; Eight CO problems; direct scheduling experiments are single-machine total weighted tardiness problem and resource-constrained project scheduling problem (RCPSP), whereas neural local search (NLS) is evaluated only on TSP/CVRP}
& Base scheduling path: SC-O; SG0; SL / adjacent routing NLS path: SC-O; SG0; SL
& RF
& \textbf{Boundary:} Retained ACO forms schedules; the neural local-search result is adjacent routing evidence only. \\
\addlinespace[1pt]

\citet{sun2023difusco} \newline \emph{Adjacent; Routing and graph CO}
& Greedy decode: FORM; SG4; MDL / greedy plus 2-opt: FORM; SG4; SH / heatmap-guided Monte Carlo tree search: SC-O; SG0; SL
& IR
& \textbf{P:} The configurations are adjacent routing evidence and do not establish scheduling performance. \\
\addlinespace[1pt]

\citet{li2026goaldiffusion} \newline \emph{Direct; Flow-shop, job-shop, and flexible-job-shop benchmarks} \newline \emph{Publication status: preprint at the cutoff}
& FORM; SG4; MDL
& IR
& \textbf{F:} Fixed decoding covers represented constraints; the direct scheduling evidence was an emerging preprint at the cutoff. \\
\addlinespace[1pt]

\citet{feng2026frontierco} \newline \emph{Direct for the FJSP benchmark; adjacent for the other CO tasks}
& Interface \textnormal{\texttt{NR}}; SG record varies by evaluated solver; authority \textnormal{\texttt{NR}}
& \texttt{NR}
& \textbf{P/A:} Only the FJSP benchmark is direct; interface and authority vary by solver and remain \texttt{NR} in the aggregate record. \\
\addlinespace[1pt]

\citet{ding2024multipolicy} \newline \emph{Direct; multi-objective multiplicity FJSP}
& FORM; SG2; MDL
& RF
& \textbf{P/A:} Produces an empirical approximation to a trade-off set; it does not certify Pareto optimality or report field use. \\
\addlinespace[1pt]

\citet{echeverria2025leveraging} \newline \emph{Direct; FJSP}
& FORM; SEL; SG2; SH
& IR
& \textbf{F:} Final validity comes from CP completion; the system is neither pure end-to-end generation nor only label generation. \\
\addlinespace[1pt]

\citet{li2025rolling} \newline \emph{Direct; long-horizon FJSP under rolling-horizon optimization}
& FORM; SG1; SH
& IR
& \textbf{F:} The learned mask simultaneously fixes assignments for one release, while the retained solver chooses remaining assignments and timing. \\
\addlinespace[1pt]

\citet{marques2025dynamic} \newline \emph{Direct; dynamic JSSP}
& SEL; SG2; MDL
& IR or RF by variant
& \textbf{P/A:} Performance is bounded by the candidate rule library; simulated dynamics do not establish deployment. \\
\addlinespace[1pt]

\citet{stockermann2025scalability} \newline \emph{Direct; semiconductor fab; Level C}
& FORM; SG2; MDL
& RF
& \textbf{A:} More than 1,000 equipment units and industrial data strengthen scale and realism, but the evaluated scheduling process remains in simulation. \\
\addlinespace[1pt]

\citet{tiacci2024simulator} \newline \emph{Enabling; dynamic FJSP simulation infrastructure}
& Interface \texttt{NR}; SG0; authority \texttt{NR}
& \texttt{NR}
& \textbf{F/A:} Provides infrastructure and a methodological critique rather than a newly trained scheduling policy; reuse improves reproducibility, not field validity. \\
\addlinespace[1pt]

\citet{yang2025goodrl} \newline \emph{Direct; dynamic workflow scheduling}
& FORM; SG2; MDL
& IR $\rightarrow$ RF $\rightarrow$ RF
& \textbf{P/A:} Explicit arrivals support a dynamic claim. \\
\addlinespace[1pt]

\citet{yedidsion2025productization} \newline \emph{Enabling; semiconductor deployment engineering}
& Interface \texttt{NR}; SG0; authority \texttt{NR}
& \texttt{NR}
& \textbf{P/A:} Describes engineering requirements but reports no named-site duration or independently evaluated field KPI. \\
\addlinespace[1pt]

\citet{zhao2024heterogeneous} \newline \emph{Direct; hybrid flow shop}
& FORM; SG2; MDL
& RF
& \textbf{A:} Results support the tested baselines and instances only; beating simple rules is not a field-wide state-of-the-art claim. \\
\addlinespace[1pt]

\citet{senoner2023distributional} \newline \emph{Direct; customized-production order fulfillment; Level B2}
& SUB; SG0; SL
& IR
& \textbf{A:} Real Aker data establish distribution-shift relevance, but the numerical/retrospective evaluation is not production deployment. \\
\addlinespace[1pt]

\citet{vaclavik2018accelerating} \newline \emph{Enabling mechanism-boundary evidence; nurse rostering and time-division-multiplexing scheduling via branch-and-price}
& \texttt{NA} (L not satisfied; retained exact search)
& IR
& \textbf{P:} The learned artifact is not evaluated beyond the solve that created it, so L is not satisfied. \\
\addlinespace[1pt]

\citet{koutecka2025rankpricing} \newline \emph{Direct; operating-room branch-and-price}
& SC-O; SG0; SL
& IR
& \textbf{A:} Learning prioritizes subproblems rather than generating the operating-room schedule; the host solver retains responsibility. \\
\addlinespace[1pt]

\citet{sarkar2025stabilized} \newline \emph{Direct; workforce scheduling}
& SEL; INIT; SG0; SL
& IR
& \textbf{Boundary:} Predictions initialize and stabilize retained column generation; pricing, master optimization, and termination remain responsible for the solution. \\
\addlinespace[1pt]

\citet{verhaeghe2024precedences} \newline \emph{Direct; resource-constrained project scheduling}
& Constraint-addition variant: FORM; SG1; SH / ordering variant: SC-O; SG0; SL
& IR
& \textbf{F/P:} Constraint addition is FORM/SG1/SH; ordering is SC-O/SG0/SL, and any certificate concerns the searched model. \\
\addlinespace[1pt]

\citet{teichteil2023resource} \newline \emph{Direct; deterministic and stochastic RCPSP}
& FORM; SG1; MDL
& IR
& \textbf{F/A:} The schedule-generation scheme provides modeled feasibility; stochastic robustness is scenario-based, not field evidence. \\
\addlinespace[1pt]

\citet{jesus2026constraint} \newline \emph{Direct; FJSP}
& FORM; SG2; MDL
& IR+RF
& \textbf{Boundary:} CP supplies training knowledge, not an inference-time exact guarantee. \\
\addlinespace[1pt]

\citet{antonov2026tardy} \newline \emph{Direct; single machine with weights, due dates, and deadlines; minimize weighted number of tardy jobs}
& PROP; SG1; SH
& IR
& \textbf{F/P/A:} The result is problem-specific and supplies neither general objective optimality nor a transfer certificate. \\
\addlinespace[1pt]

\citet{liu2026featureselection} \newline \emph{Direct; dynamic JSSP}
& SEL; SG2; MDL
& IR
& \textbf{Boundary:} Real-life-inspired events remain simulated and the fixed rule set places a ceiling on decisions. \\
\addlinespace[1pt]

\citet{uzunoglu2023learning} \newline \emph{Direct; parallel serial-batch machines}
& SEL; SG0; SL
& IR
& \textbf{P:} Scale reaches thousands of jobs, but speed and quality loss must be reported together. \\
\addlinespace[1pt]

\citet{uzunoglu2024multistart} \newline \emph{Direct; serial-batch scheduling}
& SEL; SG0; SL
& IR
& \textbf{P:} Large time reduction is accompanied by a reported solution-quality trade-off relative to full multi-start. \\
\addlinespace[1pt]

\citet{uzunoglu2025algorithmselection} \newline \emph{Direct; serial-batch scheduling}
& SEL; SG0; SL
& IR
& \textbf{Boundary:} The effect of selection must be separated from subsequent genetic-algorithm improvement; evidence is specific to this batching environment. \\
\addlinespace[1pt]

\citet{niroumandrad2024tabu} \newline \emph{Direct for the saved-model and transfer variants; physician scheduling}
& SC-C; SG0; SH
& IR
& \textbf{P/A:} Learning restricts the neighborhood, while retained tabu search selects and applies revisions. \\
\addlinespace[1pt]

\citet{wang2023simulationoptimization} \newline \emph{Direct; dynamic JSSP with uncertainty}
& SUB; SG0; SL
& IR
& \textbf{A:} The industrial case supports realism, not sustained deployment, and multi-component attribution remains difficult. \\
\addlinespace[1pt]

\citet{han2025imitation} \newline \emph{Direct; semiconductor/display-style real-time dispatching in commercial scheduling simulation}
& FORM; SG2; MDL
& IR
& \textbf{P/A:} The commercial-scheduler datasets establish imitation accuracy and latency, not teacher improvement, plant KPI gains, or deployment. \\
\addlinespace[1pt]

\citet{benda2019flowshop} \newline \emph{Direct; flow shop with alternatives, setups, and blocking}
& FORM; SG2; MDL
& IR
& \textbf{F/A:} Complex constraints broaden direct evidence; industrial motivation and simulation do not imply live use. \\
\addlinespace[1pt]

\citet{thevenin2019vns} \newline \emph{Enabling mechanism-boundary evidence; single machine with time windows and rejection}
& \texttt{NA} (L not satisfied; retained variable neighborhood search)
& RF
& \textbf{P:} Within-run trails are not evaluated in a later solve, so L is not satisfied. \\
\addlinespace[1pt]

\citet{alicastro2021rlils} \newline \emph{Enabling mechanism-boundary evidence; additive-manufacturing batch-machine scheduling}
& \texttt{NA} (L not satisfied; retained local search)
& RF
& \textbf{F:} Within-run Q-values are not evaluated beyond the solve that created them, so L is not satisfied. \\
\addlinespace[1pt]

\citet{zhang1995rljobshop} \newline \emph{Direct; job-shop-style temporal/resource scheduling, including NASA Space Shuttle Payload Processing}
& FORM; SG2; SH
& RF
& \textbf{P/A:} Historically important computational evidence; retained initialization and repair also contribute schedule content, and the experiments do not establish modern-scale performance or deployment. \\
\addlinespace[1pt]

\citet{riedmiller1999neural} \newline \emph{Direct; single- and multi-resource production scheduling}
& Single-resource: FORM; SG2; MDL; RL / multi-resource mixed-policy: FORM; SG2; SH
& RF
& \textbf{P:} Evidence is limited to small computational examples; local learning does not establish global optimality, and the mixed-resource case retains nonlearned dispatch authority. \\
\addlinespace[1pt]

\citet{aydin2000dynamic} \newline \emph{Direct; dynamic JSSP}
& SEL; SG2; MDL
& RF
& \textbf{A:} One of the early explicit dynamic studies; dynamics are simulated and the rule library limits actions. \\
\addlinespace[1pt]

\citet{gabel2012distributed} \newline \emph{Direct; deterministic/stochastic JSSP}
& FORM; SG2; MDL
& RF
& \textbf{Boundary:} Broadens policy-gradient history beyond REINFORCE/PPO, but remains computational evidence. \\
\addlinespace[1pt]

\citet{sallam2021stochastic} \newline \emph{Enabling mechanism-boundary evidence; stochastic RCPSP}
& \texttt{NA} (L not satisfied; retained metaheuristics)
& RF
& \textbf{P:} The within-run Q-table is not evaluated in a later solve, so L is not satisfied. \\
\addlinespace[1pt]

\citet{zhao2022rcpsp} \newline \emph{Direct; RCPSP}
& FORM; SG2; MDL
& RF
& \textbf{P/A:} Small-to-larger-size tests support bounded size-transfer evidence; baseline strength limits broader state-of-the-art claims. \\
\addlinespace[1pt]

\citet{cai2024rcpspdisruptions} \newline \emph{Direct; RCPSP with resource disruptions}
& FORM; SG2; MDL
& RF
& \textbf{A:} Explicit disruptions support robustness analysis; reported limitations on large instances and lack of field use remain material. \\
\addlinespace[1pt]

\citet{brammer2022permutation} \newline \emph{Direct; multi-line PFSP with demand plans}
& FORM; SG2; SH
& RF
& \textbf{P:} Short- and long-budget conclusions differ; the local-search cost belongs to inference. \\
\addlinespace[1pt]

\citet{li2022bilevel} \newline \emph{Direct; large-scale FFSP}
& FORM; SG2--SG3; MDL
& RF
& \textbf{P/A:} Public benchmarks and large industrial scenarios strengthen scale and application relevance; they do not document live schedule release or sustained use. \\
\addlinespace[1pt]

\citet{pan2023pfsp} \newline \emph{Direct; PFSP}
& FORM; SG2; SH
& RF
& \textbf{A:} Size-transfer claims must retain the exact train/test ranges and total refinement budget. \\
\addlinespace[1pt]

\citet{muller2024twostage} \newline \emph{Direct; two-stage PFSP; Level B2}
& FORM; SG2; MDL
& RF
& \textbf{A:} Five weeks of appliance-production data support retrospective realism, not documented live deployment. \\
\addlinespace[1pt]

\citet{yang2025dynamicpfsp} \newline \emph{Direct; dynamic multi-objective PFSP}
& SEL; SG2; MDL
& RF
& \textbf{P/A:} A fixed weighted objective is not a complete Pareto policy; small-to-large transfer remains benchmark-specific. \\
\addlinespace[1pt]

\citet{corsini2022permutations} \newline \emph{Direct; JSSP}
& SC-C; SG0; SH
& IR
& \textbf{P:} The learned filter restricts the neighborhood while retained tabu search chooses and applies revisions; quality and runtime concern the assembled pipeline. \\
\addlinespace[1pt]

\citet{zhang2024tbgat} \newline \emph{Direct; JSSP}
& FORM; SG3; SH
& RF
& \textbf{A:} Provides additional learned-improvement evidence across synthetic and classical instances, not deployment. \\
\addlinespace[1pt]

\citet{kaleta2026constructive} \newline \emph{Direct; FJSP}
& FORM; SG2; MDL
& RF
& \textbf{P:} Fine-tuning time is part of test-time cost; results question, rather than assume, the need for complex deep architectures. \\
\addlinespace[1pt]

\citet{yaakoubi2020machine} \newline \emph{Direct; airline crew pairing; Level B2}
& INIT; SG0; SL
& IR
& \textbf{P/A:} Six months of large-airline pairings support real-data relevance; no operational release or field KPI is reported. \\
\addlinespace[1pt]

\citet{lex2026smart} \newline \emph{Direct; operating-room scheduling; Level B2}
& SUB; SG0; SL
& IR
& \textbf{Boundary:} The retrospective cohort contains 15,267 cases; weekly schedules were generated and simulated rather than released, and benefits depend on both prediction and the retained integer linear program. \\
\addlinespace[1pt]

\citet{wallner2025hierarchical} \newline \emph{Direct; operating-room fairness; Level B2}
& SUB; FORM; SG2; MDL
& IR+RF
& \textbf{Boundary:} Seventy-eight weeks of real schedules support realism; proposed allocations were not documented as adopted. \\
\addlinespace[1pt]

\citet{chorfi2026analytics} \newline \emph{Enabling evaluation evidence; nurse scheduling under absence uncertainty; Level B2}
& No release-path learned interface; SG and authority \texttt{NA}
& IR
& \textbf{A:} Learning generates evaluation scenarios rather than release-path schedule content; architecture is \texttt{NA}. \\
\addlinespace[1pt]

\citet{zhang2026digitaltwin} \newline \emph{Direct; dynamic framework tested on a hybrid-flow-shop-like water-heater welding line; Level B1}
& FORM; SG2; MDL
& RF
& \textbf{P/A/R:} The physical workshop experiment supports B1 pilot evidence; sustained release, operator use, fallback, and before/after plant KPIs are not documented. \\
\addlinespace[1pt]

\citet{mao2019decima} \newline \emph{Direct; data-processing DAG scheduling; Level B1}
& FORM; SG2; MDL
& RF
& \textbf{A:} A 25-node testbed and production-trace replay are stronger than simulation, but not production deployment. \\
\addlinespace[1pt]

\citet{jajoo2022slearn} \newline \emph{Direct; cluster/DAG scheduling; Level B1}
& SUB; SG0; SL
& IR
& \textbf{A:} Three production traces and an Azure testbed support external realism; results concern the assembled estimator--scheduler system, and no daily release is claimed. \\
\addlinespace[1pt]

\citet{delimitrou2013paragon} \newline \emph{Direct; heterogeneous datacenter scheduling; Level B1}
& SUB; SG0; SL
& IR
& \textbf{A:} Controlled cloud and physical-system experiments support B1 testbed evidence; sustained tenant-facing production deployment is not documented. \\
\addlinespace[1pt]

\citet{narayanan2020gavel} \newline \emph{Direct; heterogeneous graphics processing unit (GPU) scheduling; Level B1}
& SUB; SG0; SL
& IR
& \textbf{A:} Physical-cluster experiments support B1 testbed evidence; estimator value depends on retained optimization and scheduling, and production release is not documented. \\
\addlinespace[1pt]

\citet{heik2024industrial} \newline \emph{Direct; dynamic manufacturing; Level C}
& FORM; SG2; MDL
& RF
& \textbf{A:} The physical test-bed motivates the twin, but the reported scheduling experiments are digital rather than live. \\
\addlinespace[1pt]

\citet{cotopalacio2022qlearning} \newline \emph{Direct; real-company human--robot FJSP; Level C}
& FORM; SG2; MDL
& RF
& \textbf{A:} The company scenario is concrete, but physical execution and sustained adoption are not reported. \\
\addlinespace[1pt]

\citet{han2024wetclean} \newline \emph{Direct; semiconductor wet-clean station; Level C}
& SUB; SG2; MDL
& IR
& \textbf{A:} The study therefore supports Level C industrial simulation, not retrospective operational-data evaluation or field release. \\
\addlinespace[1pt]

\end{longtable}
\endgroup

Because the analytical library is purposive, neither the registry composition
nor internal code counts estimate the prevalence of methods, architectures, or
assurance controls.
\section{Claim-Level Technical Evidence Audit}
\label{ec:technical_evidence_audit}

This section supplies the technical claim--evidence records behind Section~5
of the main paper. The unit is a claim about a reported system configuration,
not a model result in isolation. Each row identifies the taxonomy boundary
that determines the object of the claim, the sources that bear on it, and the
maximum supported conclusion and excluded inference. Different objectives,
instance generators, release paths, and computing budgets are not pooled into
a common effect size. Operational maturity is assessed separately in EC.5.

\subsection{Evidence Matrix for the Main Technical Conclusions}
\label{ec:technical_claim_matrix}

\begingroup
\setlength{\tabcolsep}{2pt}
\footnotesize
\begin{longtable}{L{2.30cm}L{4.55cm}L{4.05cm}L{4.55cm}}
\caption{Claim-level technical evidence matrix}
\label{tab:ec_technical_claim_matrix}\\
\toprule
Technical conclusion
& Taxonomy and system boundary
& Evidence sources
& Maximum supported conclusion and excluded inference \\
\midrule
\endfirsthead
\multicolumn{4}{l}{\small\itshape Table \thetable\ continued}\\
\toprule
Technical conclusion
& Taxonomy and system boundary
& Evidence sources
& Maximum supported conclusion and excluded inference \\
\midrule
\endhead
\bottomrule
\endfoot

System-equivalence boundary in quality--computation-time comparisons
& The compared object includes the released decision, information and
objective, complete normal path, release requirement, and decision budget.
Direct formation still includes decoding and checking; a long solver run used
as a reference is not a time-matched release path.
& Learning to Dispatch and the heterogeneous-graph method compare with stated
dispatching rules and report OR-Tools or CPLEX reference solutions under limits
of up to 3,600 and 1,800 seconds, respectively
\citep{zhang2020dispatch,zhao2024heterogeneous}. Imitation learning compares
with its NEH teacher and other reported methods \citep{li2024imitation}.
& The records support gains over the listed rules and a size-dependent
quality--computation-time trade-off against the reported teacher. They do not
establish a general ranking over strong optimization, metaheuristic, or learned
alternatives under a common release contract. \\

Reported work and cost displacement
& Work is assigned to capability acquisition, per-release formation or search,
retained optimization, checking, and adaptation; integration, monitoring,
human intervention, and maintenance are recorded as unreported when absent.
& Teacher generation and training are reported for imitation learning
\citep{li2024imitation}; SLIM evaluates model-generated candidates during
self-labeling \citep{corsini2024selflabeling}; MISTAR evaluates repeated
improvement moves \citep{wang2025mistar}; and rolling-horizon learning calls a
model and retained subsolver at successive releases \citep{li2025rolling}.
& The sources show that reported computation can move among acquisition,
formation, improvement, and retained optimization rather than disappear. They
do not provide a common estimate of total decision cost, break-even reuse, or
lifecycle value. \\

Guarantee localization
& Technical authority identifies the enforcing component and the decision
space it can control. Order-preserving guidance leaves the declared search
space available; a learned machine order or fixed assignments restrict what a
downstream optimizer can revise.
& Learned ordering within complete CP search is SG0 and solver-led
\citep{sun2024enhancing}. Parallel prediction with LP completion and
rolling-horizon learned fixing are SG1 paths with shared technical authority
\citep{kotary2022fast,li2025rolling}.
& Certified CP termination can support the declared solver claim. Downstream
feasibility, bound, or optimality statements in the shared-authority paths are
confined to the formulation and residual decision space; solver presence does
not extend them to choices already fixed by learning. \\

Same-artifact within-family transfer
& A transfer record identifies the retained learned artifact, scheduling
family, changed size, decoder, comparator, adaptation budget, and post-change
feasibility, quality, and decision time. Execution alone establishes only
computational reach.
& Trained policies are evaluated on larger JSSPs
\citep{zhang2020dispatch,park2021learning}; a policy trained on 20-job PFSP
instances is tested up to 1,000 jobs \citep{li2024imitation}; and separately
trained RESCHED configurations are tested at larger sizes within their
respective families \citep{xiao2026resched}.
& These studies support bounded same-artifact size transfer within the
reported families and ranges. They do not establish unrestricted size
invariance or structural transfer across formulations or tasks. \\

Cross-problem transfer identity
& Same-artifact transfer, adaptation-based transfer, architecture reuse after
separate training, and strategy diversity within one task are distinct
records.
& RESCHED trains separate FJSP, JSSP, and FFSP configurations under a common
architecture \citep{xiao2026resched}. GOAL fine-tunes a multitask backbone
trained partly on JSSP for open-shop scheduling \citep{drakulic2025goal}.
PolyNet generates strategy diversity within FFSP \citep{hottung2025polynet}.
& GOAL supplies a bounded adaptation-based cross-problem record. RESCHED
supports architecture reuse and within-family evaluation; PolyNet supports
within-task diversity. None establishes general zero-shot one-policy transfer
across scheduling formulations or tasks. \\

Dynamic operating regime
& The supported regime is defined by the arrival process, information at each
decision, replanning epoch, simulator, update phase, and complete release-path
computation. Repeated fixed-policy use is response within a regime; updating
the policy is adaptation.
& The studies evaluate repeated rule selection under new-job insertion
\citep{luo2020dynamic}, imitation and PPO under workflow arrivals
\citep{yang2025goodrl}, and learned fixing with retained rolling-horizon
optimization \citep{li2025rolling}.
& The systems improve reported outcomes under the declared event and update
processes. They do not establish adaptation to an unmodeled event process,
long-run stability under persistent change, or deadline-feasible release
outside the tested regime. \\

Stochastic and information-reveal regime
& The uncertainty distribution, information revealed before each action,
scenario and evaluation budgets, risk functional, decoder, and comparator
limit define the claim.
& Expected makespan and a stated value-at-risk measure are evaluated under
sampled processing-time scenarios \citep{smit2025stochastic}. Dispatching
under a declared structural-uncertainty and reveal process optimizes expected
performance and reports conditional value at risk as an evaluation measure
\citep{zhang2026uncertainty}.
& The records support expected or risk-related performance under the reported
uncertainty and reveal regimes. They do not establish robustness to a different
distribution, information process, scenario budget, or risk criterion. \\

Preference and objective regime
& The objective set, scalarization or preference input, tested weight domain,
policy collection, decoder, comparator budget, and quality indicator delimit
the supported regime. Pareto coverage and preference generalization are
empirical claims.
& The evidence covers dynamic rule selection for alternative scalar objectives
\citep{luo2021multiobjective}, multiple policies for stated weight settings
\citep{ding2024multipolicy}, and a preference-conditioned policy evaluated on
untrained weights within a specified domain \citep{smit2026multiobjective}.
& The studies support trade-off generation and bounded preference
generalization for the tested objectives and domains. They do not certify
Pareto completeness, validate stakeholder preferences, or establish transfer
after a change in constraints, resources, or scheduling formulation. \\

\end{longtable}
\normalsize
\endgroup

\section{Operational Maturity and Source-Adjudicated Application Episodes}
\label{ec:operational_evidence}

An application episode is a reported system--organization--setting
combination coded from the full text. \emph{Source-adjudicated} means that the
record rests on explicit source reporting, not that reported use or outcomes
were independently verified. Provenance is marked \textbf{F} when the episode
was recovered by the structured direct-study route and \textbf{P} when it
entered through purposive supplementation; provenance does not determine the
operational-evidence level.

\subsection{Operational-Maturity Rules}
\label{ec:operational_maturity}

Operational maturity records the closest documented contact between the
evaluated system and operating decisions. It is coded independently of
inference architecture, normal-path technical authority, automated release,
and organizational decision rights. Company participation, operational data,
an industrial problem description, or industry-calibrated simulation alone
does not establish Level A or Level B1 evidence.

\begin{table}[!htbp]
\centering
\caption{Operational-evidence ladder}
\label{tab:ec_evidence_ladder}
\small
\begin{tabularx}{\textwidth}{L{1.1cm}L{6.0cm}Y}
\toprule
Level & Minimum documented evidence & Interpretation boundary \\
\midrule
A
& Routine or sustained use in an ongoing operating process, with the reported
release or human-approval arrangement identified.
& Establishes operating use for the reported period, not autonomy, causal
impact, continued use after publication, or lifecycle value. \\
B1
& Bounded physical or end-to-end execution, physical pilot, shadow operation,
user-acceptance exercise, or controlled testbed under stated conditions.
& Establishes bounded execution or organizational contact, not routine or
sustained use. \\
B2
& Retrospective evaluation on operational data without physical execution or
release of the evaluated schedules.
& Establishes relevance to observed operations, not prospective use. \\
C
& Industry-calibrated model or simulation representing domain-specific
equipment, constraints, and events.
& Establishes industrial realism within the model, not field use. \\
D
& Synthetic instances or stylized simulation.
& Establishes behavior within the reported experimental design. \\
\bottomrule
\end{tabularx}
\end{table}

\subsection{Unified Application-Episode Register}
\label{ec:formal_route_episodes}
\label{ec:supplemental_episodes}

The register reports Level A and B1 episodes because these records establish
operating use or bounded end-to-end execution. Rights are specification (S),
release (R), intervention (I), and lifecycle (L); a holder is named only when
the source connects the role to the right, and otherwise is recorded as
\texttt{NR}.

\begingroup
\setlength{\tabcolsep}{2.5pt}
\renewcommand{\arraystretch}{1.06}
\scriptsize
\begin{longtable}{L{2.65cm}L{4.15cm}L{4.35cm}L{3.85cm}}
\caption{Source-adjudicated operational-maturity episode register}
\label{tab:ec_operational_episodes}
\label{tab:ec_formal_episodes}
\label{tab:ec_supplemental_episodes}\\
\toprule
Episode, provenance, and level
& Configuration, technical authority, and operating evidence
& Assurance, recovery, and reported organizational rights
& Maximum claim and source boundary \\
\midrule
\endfirsthead
\multicolumn{4}{l}{\small\itshape Table \thetable\ continued}\\
\toprule
Episode, provenance, and level
& Configuration, technical authority, and operating evidence
& Assurance, recovery, and reported organizational rights
& Maximum claim and source boundary \\
\midrule
\endhead
\bottomrule
\endfoot

\multicolumn{4}{l}{\textbf{Panel A. Level A operating-use episodes}} \\
\addlinespace[2pt]

\citet{kim2024photolithography}\newline
\emph{AMOLED photolithography; \textbf{F}; learning-augmented; Level A}
& When the learned candidate is selected, learned constructive decisions are
binding and the path is model-led; selection of the retained rule-based
candidate yields a solver-led path. Both candidates ran in actual fabrication
for one year. The published comparison covers four consecutive months at about
480 decisions per day.
& Action filters remove unavailable actions, and a fixed KPI selector chooses
between the normal-path candidates before dispatch; this is not documented
failure-triggered recovery. Rights: S---\texttt{NR}; R---automated selector,
formal holder \texttt{NR}; I---\texttt{NR}; L---\texttt{NR}.
& Establishes sustained operating use under the reported dual-candidate path.
The source reports improvements in step-target matching, setup, and assignment
measures, but not complete intervention, failure, maintenance, rollback, or
randomized-counterfactual records. \\

\citet{liang2022lenovo}\newline
\emph{Lenovo/LCFC laptop assembly; \textbf{P}; learning-augmented; Level A}
& A DRL model forms assignment and sequence content; planners review the
schedule and may lock manual changes before remaining work is rescheduled. The
platform was rolled out to four plants and 43 lines. Reported planning time
fell from about six hours to 30 minutes, backlog by 20\%, and fulfillment
improved by 23\%.
& Feasibility and preference controls are embedded in the platform; planner
review and modification are reported, but failure-triggered recovery is not.
Rights: S---planners set preferences; R---planner review; I---planners may
modify and lock work; L---\texttt{NR}.
& Establishes operating use with documented planner rights. The before--after
account does not isolate causal effects, and complete development,
maintenance, failure, override, and lifecycle-cost records are not reported. \\

\citet{hong2024naphtha}\newline
\emph{LG Chem Daesan naphtha-cracking center; \textbf{P};
learning-augmented; Level A}
& A multi-agent policy, simulator, and beam search form candidate schedules;
plant staff select the schedule used in operation. The computational candidate
path is model-led, while staff hold the organizational release right. The
source states that the online scheduling service is in operational use.
& Staff upload inputs, inspect candidates, and download the selected schedule;
no separate failure-triggered fallback is reported. Rights: S---staff provide
operating inputs; R---plant staff; I---staff review, further authority
\texttt{NR}; L---\texttt{NR}.
& Establishes operating use, not a field performance comparison. Reported
profit and constraint comparisons are simulator-based; field duration,
overrides, recovery rates, and a field counterfactual are not reported. \\

\citet{sood2024fab}\newline
\emph{Micron semiconductor fabs; \textbf{P}; learning-augmented; Level A}
& A staged configuration uses shadow recommendations and shift-level engineer
decisions; a later configuration is reported as an autonomous, model-led path.
The report describes production deployment, user-acceptance testing,
monitoring, retraining triggers, and rollout.
& Manufacturing engineers set objective weights and initially decide whether
to follow recommendations; release is reported as automated in the later
phase. Recovery procedures are not specified. Rights: S---manufacturing
engineers; R---engineers in the staged phase and automated later;
I---\texttt{NR}; L---\texttt{NR}.
& Establishes staged and reported autonomous production use. Readable effect
magnitudes, duration of autonomous use, failure and override rates, recovery
performance, and lifecycle cost are not supplied. \\

\citet{konstantelos2022fabwide}\newline
\emph{Seagate Springtown fab; \textbf{P}; solver-led industrial reference;
Level A}
& Hierarchical fab-wide and toolset schedulers were trialed from March to May
2022 and reported as operating continuously from June 2022. Fab-state data
arrive every five minutes; local schedules take about five minutes and the
fab-wide run generally takes about seven minutes.
& Managers progressively expanded tool coverage and guidance strength as
operators validated decisions. Rights: S---managers configure coverage and
guidance strength; R---formal holder \texttt{NR}; I---operators and managers
participate in validation, exact allocation \texttt{NR}; L---\texttt{NR}.
& Establishes a solver-led Level A reference, not a matched comparator for the
four learning-augmented Level A episodes. Use beyond the reporting period and
complete failure, rollback, counterfactual, and lifecycle-cost records are not
established. \\

\addlinespace[3pt]
\multicolumn{4}{l}{\textbf{Panel B. Level B1 bounded experiments and testbeds}} \\
\addlinespace[2pt]

\citet{zhang2026digitaltwin}\newline
\emph{Water-heater welding line; \textbf{F}; learning-augmented; Level B1}
& A learned constructive policy forms binding decisions on a model-led path.
After virtual training, the system read and wrote physical physical programmable logic controller data and issued
commands on the line. One experiment lasted about 13 minutes and handled 23
requests.
& Implemented decision rules constrain the physical commands; a separate
checker or fallback is not documented. Rights: S/R/I/L---\texttt{NR}.
& Establishes bounded physical release, not routine use. Long-run human roles,
failures, recovery, maintenance, and a matched operating counterfactual are not
reported. \\

\citet{narayanan2020gavel}\newline
\emph{Heterogeneous GPU cluster; \textbf{F}; learning-augmented; Level B1}
& Online profiling and matrix completion estimate job--accelerator throughput
for retained optimization and round-based scheduling. The complete physical
research-cluster testbed path is solver-led.
& Retained optimization and scheduling form the release-level decisions;
organizational release and recovery procedures are not documented. Rights:
S/R/I/L---\texttt{NR}.
& Establishes end-to-end testbed execution, not routine production release,
sustained recovery performance, or lifecycle maintenance. \\

\citet{mao2019decima}\newline
\emph{Data-processing cluster; \textbf{P}; learning-augmented; Level B1}
& A learned constructive policy forms the schedule on a model-led path. The
scheduler was executed on a 25-node cluster testbed with workloads generated
from production traces.
& The policy forms the testbed schedule; organizational release, independent
checking, and recovery arrangements are not documented. Rights:
S/R/I/L---\texttt{NR}.
& Establishes physical testbed execution and trace realism, not routine use,
long-run failure behavior, or maintenance. \\

\citet{jajoo2022slearn}\newline
\emph{Azure cluster testbed; \textbf{P}; learning-augmented; Level B1}
& Task sampling estimates runtime properties for a retained scheduler evaluated
with three production traces on an Azure testbed. Because the retained
scheduler forms release-level content, the path is solver-led.
& The retained scheduler forms the evaluated decisions; organizational
release and fallback are not documented. Rights: S/R/I/L---\texttt{NR}.
& Establishes the complete estimator--scheduler testbed path, not routine
service use, sustained governance, or lifecycle outcomes. \\

\citet{delimitrou2013paragon}\newline
\emph{Heterogeneous datacenter testbed; \textbf{P}; learning-augmented;
Level B1}
& Brief profiling and collaborative filtering estimate workload heterogeneity
and interference for a retained greedy placement scheduler. The solver-led
path was evaluated in controlled cloud and physical-system experiments.
& The retained scheduler determines schedule content; sustained release and
failure-triggered recovery are not documented. Rights:
S/R/I/L---\texttt{NR}.
& Establishes bounded testbed execution, not sustained tenant-facing
production use, failure and recovery rates, or maintenance effort. \\

\end{longtable}
\normalsize
\endgroup

The register contains four learning-augmented Level A operating-use episodes,
one solver-led Level A industrial reference case, and five learning-augmented
Level B1 experiments or testbeds. It establishes operational possibility under
the reported conditions, not deployment prevalence, causal superiority, or
long-run lifecycle value. Operational maturity does not by itself identify
technical authority, automated release, recovery capacity, or the holders of
specification, release, intervention, and lifecycle rights.
\section{Evidence-to-Proposition Trace}
\label{ec:proposition_trace}

Section~6 of the main paper states the four propositions and their
operations-management mechanisms. The review establishes feasible
configurations, bounded outcomes, and unresolved comparisons, but the reviewed
scheduling studies do not directly estimate the proposed relationships. This
section therefore records the source basis and remaining gap for each
proposition and specifies a direct test.

Tests should compare complete systems with a prespecified best feasible
alternative over a common accounting period. \emph{Relative lifecycle value}
is the difference in operating benefit from released decisions after deducting
capability-acquisition, integration, decision-production, assurance, failure,
recovery, updating, and maintenance costs. Safety, legal, and continuity
requirements remain hard constraints. The evaluation contract should fix or
measure the problem and event distribution, decision-time information,
objectives and constraints, release requirement and deadline, computing and
human resources, comparator tuning, adaptation budget, and accounting horizon.
For each directional estimand below, a sign reversal or an upper confidence
bound below a prespecified smallest meaningful positive effect counts against
the proposition; imprecision alone does not. Failure of a stated boundary
condition is not contrary evidence.

\begingroup
\setlength{\tabcolsep}{3pt}
\renewcommand{\arraystretch}{1.05}
\scriptsize
\begin{longtable}{L{1.00cm}L{4.25cm}L{6.15cm}L{3.80cm}}
\caption{Reviewed gaps, constructs, and direct tests for the four propositions}
\label{tab:ec_prop_tests}\\
\toprule
Prop.
& Reviewed basis and unresolved relationship
& Constructs, design, and estimand
& Evidence inconsistent with the proposition \\
\midrule
\endfirsthead
\multicolumn{4}{l}{\small\itshape Table \thetable\ continued}\\
\toprule
Prop.
& Reviewed basis and unresolved relationship
& Constructs, design, and estimand
& Evidence inconsistent with the proposition \\
\midrule
\endhead
\bottomrule
\endfoot

P1
& The Lenovo, photolithography, and Seagate records document repeated use and
supporting data, validation, or integration assets
\citep{liang2022lenovo,kim2024photolithography,
konstantelos2022fabwide}. They do not provide common lifecycle accounting or
estimate how value changes with effective reuse.
& \emph{Effective reuse volume} is the number of decisions covered during the
joint operating range and useful life of the learned capability and its
supporting assets. Measure full acquisition and integration cost, recurring
complete-path benefit and cost, and updating and maintenance. Using
longitudinal, staggered-adoption, or matched multi-site data, estimate the
change in relative lifecycle value as effective reuse volume grows,
conditional on a positive recurring advantage over the best feasible
alternative.
& With a positive recurring advantage and the stated coverage maintained,
relative lifecycle value decreases materially as effective reuse volume grows.
\\

P2
& Direct formation, learned warm starts, and learned fixing leave different
amounts of decision space for retained optimization
\citep{omalley2023hybrid,li2025rolling}. Existing comparisons do not jointly
vary the ex ante bottleneck, binding learned authority, retained revisable
decision space, and form of change.
& Measure the scope of choices fixed by learning, the decision space retained
for explicit revision, and the ex ante share of critical-path response time
spent on schedule formation. Distinguish \emph{state change}, which varies
information within a maintained formulation, from \emph{structural change},
which alters objectives, constraints, routing, resource relations, or another
element of the formulation. In a factorial simulation, controlled testbed, or
matched field design, estimate two moderation effects: whether a formation
bottleneck and stable formulation increase the value of expanding binding
learned authority, and whether structural change increases the value of
retained revisable authority. Hold the release contract and adaptation budget
common.
& Either predicted moderation reverses materially under its stated boundary:
expanding binding learned authority loses value at a verified formation
bottleneck under a stable formulation, or retained revisable authority loses
value as structural change increases. \\

P3
& Photolithography reports action filters and parallel candidates, the naphtha
episode staff release, and Seagate staged validation
\citep{kim2024photolithography,hong2024naphtha,
konstantelos2022fabwide}. These arrangements do not estimate the joint effect
of release scope, error consequences, assurance coverage, and recovery
capacity.
& Record the decisions eligible for automated release, error incidence and
severity, coverage and error rates for feasibility, performance, and
applicability assurance, and the cost and tail latency of fallback, rechecking,
and recovery. In simulation, shadow operation, or bounded release, vary release
scope, error consequences, assurance coverage, recovery capacity, and the
deadline; estimate their interactions in relative lifecycle value and report
average and tail outcomes. Recovery is deadline-feasible only when an available
process can produce and recheck an acceptable alternative within the deadline.
& One or more predicted moderation effects reverses: wider release loses
relative lifecycle value as error consequences fall, or as effective assurance
coverage and demonstrably deadline-feasible recovery capacity increase, after
their added costs are included. A missed deadline when recovery capacity is not
deadline-feasible is a failed boundary condition, not contrary evidence. \\

P4
& Field episodes assign problem inputs, schedule approval, intervention, and
system change to different roles
\citep{liang2022lenovo,kim2024photolithography,
konstantelos2022fabwide}. No reviewed study holds the technical pipeline fixed
while reallocating specification, release, intervention, and lifecycle rights;
unreported rights are coded \texttt{NR}.
& For each right, record its formal holder, effective control, relevant
information and expertise, accountability, and a prespecified alignment score;
also measure coordination latency, escalation, intervention, and rollback time
and success. Reallocate rights in a randomized or staggered design, or compare
matched sites while holding the technical pipeline and release contract fixed.
Estimate the change in relative lifecycle value associated with improved
alignment and its interactions with timely coordination and effective
rollback.
& With the technical pipeline fixed, improved alignment reduces relative
lifecycle value even though required coordination remains timely and effective
rollback is available. Coordination or rollback failure outside those boundary
conditions is not contrary evidence. \\

\end{longtable}
\normalsize
\endgroup

\endgroup

\clearpage
\ifpomstemplate
  \bibliographystyle{pomsref}
\else
  \bibliographystyle{plainnat}
\fi
\bibliography{poms_revised_references}

\end{document}